\documentclass[12pt]{article}
\usepackage{amsmath,amssymb,fullpage}
\usepackage[applemac]{inputenc}

\newenvironment{myenumerate}{%

\begin{enumerate}}{\end{enumerate}}

\newcommand{\dproof}{\noindent {Proof.} \quad}
\newcommand{\fproof}{\hfill $\square$ \bigskip}

\newtheorem{example}{Example}[section]
\newtheorem{theorem}[example]{Theorem}

\newtheorem{definition}[example]{Definition}
\newtheorem{lemma}[example]{Lemma}

\newtheorem{problem}[example]{Problem}
\newtheorem{remark}[example]{Remark}
\newtheorem{corollary}[example]{Corollary}

\numberwithin{equation}{section}

\def\RB{\mathbb{R}}

\def\FB{\mathbb{F}}

\def\UB{\mathbb{U}}

\def\GB{\mathbb {G}}

\def\FC{\mathcal{F}}
\def\PC{\mathcal{P}}
\def\EC{\mathcal{E}}
\def\AC{\mathcal{A}}

\def\CC{\mathcal C}
\def\RC{\mathcal R}
\def\HC{\mathcal H}
\def\R{{\bf R}}

\def\WC{\mathcal W}

\def\LC{\mathcal L}
\def\TC{\mathcal T}
\def\Scal{\mathcal S}
\def\GC{\mathcal G}

\def\ess{\mathop{ess \;sup}}

\def\1B{\text{1\!\!I}}

\def\tY{\tilde{Y}}
\def\tN{\tilde{N}}

\def\tc{\tilde{c}}
\def\tu{\tilde{u}}

\def\tZ{\tilde{Z}}
\def\tK{\tilde{K}}

\def\hX{\hat{X}}
\def\hZ{\hat{Z}}

\def\hY{\hat{Y}}
\def\hu{\hat{u}}

\def\hb{\hat{b}}

\def\hg{\hat{g}}
\def\hc{\hat{c}}
\def\hQ{\hat{Q}}
\def\hK{\hat{K}}

\def\hH{\hat{H}}
\def\hf{\hat{f}}

\def\hp{\hat{p}}
\def\hth{\hat{\theta}}
\def\hq{\hat{q}}
\def\hr{\hat{r}}
\def\hw{\hat{w}}
\def\hR{\hat{R}}

\def\hla{\hat{\lambda}}
\def\hsi{\hat{\sigma}}
\def\hpi{\hat{\pi}}
\def\hga{\hat{\gamma}}

\def\hel{\hat{\ell}}

\begin{document}

\title{Risk minimization in financial markets modeled by It\^o-L\'evy processes}

\author{
Bernt \O ksendal$^{1}$ \and Agn\`es Sulem$^{2}$}

\date{30 March 2014}

\footnotetext[1]{Department of Mathematics, University of Oslo, P.O. Box 1053 Blindern, N--0316 Oslo, Norway, and 
   Norwegian School of Economics (NHH), Helleveien 30, N--5045 Bergen, Norway.  
email: {\tt oksendal@math.uio.no}. \\
The research leading to these results has received funding from the European Research Council under the European Community's Seventh Framework Programme (FP7/2007-2013) / ERC grant agreement no [228087].}
\footnotetext[2]{ INRIA Paris-Rocquencourt, Domaine de Voluceau, BP 105, Le Chesnay Cedex, 78153, France;   \\
Universit\'e Paris-Est, F-77455 Marne-la-Vallée, France, and Department of Mathematics, University of Oslo, email: {\tt agnes.sulem@inria.fr}}
\maketitle

\paragraph{MSC(2010):} 60H10, 60H20, 60J75, 93E20, 91G80, 91G10, 91A23, 91B70, 91B30

\paragraph{Keywords:} Convex risk measure, risk minimization, recursive utility, utility optimization, It\^o-L\'evy process, backward stochastic differential equation, the maximum principle for stochastic control of FBSDE's, stochastic differential game, HJBI equation.

\begin{abstract}

This paper is mainly a survey of recent research developments regarding methods for risk minimization in financial markets modeled by It\^o-L\'evy processes, but it also contains some new results on the underlying stochastic maximum principle. 

The concept of a convex risk measure is introduced, and two representations of such measures are given, namely : 
(i) the dual representation and 
(ii) the representation by means of backward stochastic differential equations (BSDEs) with jumps.
Depending on the representation, the corresponding risk minimal portfolio problem is studied, either in the context of
 stochastic differential games
or
 optimal control of forward-backward SDEs.

The related concept of recursive utility is also introduced, and corresponding recursive utility maximization problems are studied.

In either case the maximum principle for optimal stochastic control plays a crucial role, and in the paper we prove a version of this principle which is stronger than what was previously known.

The theory is illustrated by examples, showing explicitly the risk minimizing portfolio in some cases.

\end{abstract}

\section*{Introduction}\label{sec1}

In the recent years there has been an increased focus on the concepts of risk and methods for risk minimization in finance. The purpose of this paper is to give a brief survey of this topic, and its relation to backward stochastic differential equations (BSDEs), stochastic control of forward-backward stochastic differential equations (FBSDEs) and stochastic differential games, all within the context of financial markets modeled by It\^o-L\'evy processes.

Here is an outline of the paper:
\tableofcontents

\section{BSDEs, convex risk measures and recursive utilities}\label{sec2}

In this chapter we give an introduction to backward stochastic differential equations (BSDEs) with jumps, and we relate them to the concepts of {\it recursive utilities} and {\it convex risk measures}. This section, in particular the results on BSDEs with jumps and on dynamic risk measures  is based on the paper \cite{QS}. For a similar introduction in the Brownian motion  case, we refer the reader to the survey paper on BSDEs by M.C. Quenez \cite{Q}.

From now on we let  $B(t)$ and $\tN(dt,d\zeta):= N(dt,d\zeta) -\nu(d\zeta)dt$ denote a Brownian motion and an independent  compensated Poisson random measure, respectively, on a filtered probability space $(\Omega,   \FC, \FB:=\{\FC_t\}_{0 \leq t \leq T}, P)$ satisfying the usual conditions,  $P$ is a reference probability measure and $\nu$ is the Lévy measure of $N$.

\subsection{Examples}\label{sec2.1}
We first give some examples where BSDEs appear. For simplicity we do not include jumps in these examples. The more general versions with jumps will be discussed in the subsequent sections.

\begin{example}(Replicating portfolio)\label{ex2.1} \rm
Consider a financial market with one risk free and one risky investment possibility, with prices $S_0(t), S_1(t)$ per unit given by, respectively
\begin{equation}\label{eq2.1}
\begin{cases}
dS_0(t) = S_0(t) r(t) dt \; ; \; S_0(0) = 1 \\
dS_1(t) = S_1(t) [\mu(t) dt + \sigma(t) dB(t)] \; ; \; S_1(0) > 0.
\end{cases}
\end{equation}
Let $\pi(t)$ be a self-financing portfolio, representing the fraction of the total wealth $Y(t) = Y_\pi(t)$ invested in the risky asset at time $t$. The corresponding wealth process $Y(t)$ is given by
\begin{align}\label{eq2.2}
dY(t) & = \frac{(1 - \pi(t))Y(t)}{S_0(t)} dS_0(t) + \frac{\pi(t) Y(t)}{S_1(t)} dS_1(t) \nonumber \\
& = Y(t) [\{ (1 - \pi(t)) r(t) + \pi(t) \mu(t)\}dt + \pi(t) \sigma(t) dB(t)].
\end{align}
Let $F \in L^2(\FC_T,P)$ be a given $T$-claim. We want to find $Y(0) = y > 0$ and $\pi(t)$ such that
\begin{equation}\label{eq2.3}
Y(T) = F \; \text{ a.s. }
\end{equation}
Put
\begin{equation}\label{eq2.4}
Z(t) = Y(t) \pi(t) \sigma(t).
\end{equation}
Then
\begin{equation}\label{eq2.5}
\pi(t) = \frac{Z(t)}{Y(t) \sigma(t)}
\end{equation}
and \eqref{eq2.2} becomes
\begin{equation}\label{eq2.6}
dY(t) = \left\{ r(t) Y(t) + \frac{Z(t)}{\sigma(t)} (\mu(t) - r(t))\right\} dt + Z(t) dB(t) \; ; \; 0 \leq t \leq T.
\end{equation}

The pair \eqref{eq2.6}, \eqref{eq2.3} of equations is an example of a (linear) BSDE in the pair $(Y(t), Z(t))$ of unknown processes. If we can solve this equation for $Y(t), Z(t)$, then the replicating portfolio $\pi(t)$ is given by \eqref{eq2.5}.

Note that, in contrast to ordinary SDEs in one unknown process, this equation has {\it two} unknown processes and the {\it terminal} value $Y(T)$ of $Y$ is given, not the initial value.

More generally, let
$$g(t,y,z,\omega) : [0,T] \times \RB \times \RB \times \Omega \rightarrow \RB$$
be an $\FC_t$-adapted stochastic process in $(t,\omega)$ for each $y,z$. Then the equation
\begin{equation}\label{eq2.7}
\begin{cases}
dY(t)  = - g(t, Y(t), Z(t), \omega)dt + Z(t) dB(t) \; ; \; 0 \leq t \leq T \\
Y(T) = F \text{ a.s. }
\end{cases}
\end{equation}
is a BSDE in the unknown $\FC_t$-adapted processes $(Y(t), Z(t))$ (driven by Brownian motion). (See Section \ref{sec2.2} for a more comprehensive presentation).

For simplicity of notation we suppress $\omega$ in the following.
\end{example}

\begin{example}[Recursive utility]\label{ex2.2} 
(Duffie \& Epstein (1992), Epstein \& Zin (1989), Kreps \& Porteus (1978)).
\rm
Let $g(t,y,c)$ be an $\FC_t$-adapted process. Assume that $c \rightarrow g(t,y,c)$ is concave for all $t,y$, and let $F$ be a claim in $L^2({\cal F}_T)$. Then the {\it recursive utility} process of a given consumption process $c(\cdot) \geq 0$ is defined as the solution $Y(t) = Y_g(t)$ of the equation
\begin{equation}\label{eq2.8}
Y(t) = E \left[ \int_t^T g(s,Y(s), c(s))ds \mid \FC_t \right] \; ; \; 0 \leq t \leq T.
\end{equation}
In Section \ref{sec2.2} we shall see that \eqref{eq2.8} is equivalent to the following BSDE in $(Y_g, Z_g)$:
\begin{equation}\label{eq2.9}
\begin{cases}
dY_g(t) = - g(t, Y(t), c(t)) dt + Z_g(t) dB(t) \\
Y_g(T) = 0.
\end{cases}
\end{equation}

In particular, the (total) {\it recursive utility} $U(c)$ of a given consumption process $c(\cdot)$ is defined as
\begin{equation}\label{eq2.10}
U(c) := Y_g(0).
\end{equation}
\end{example}

\begin{example}[Convex risk measures]\label{ex2.3}
(F\"ollmer \& Schied (2002), Frittelli \& Rosazza-Gianin (2002))
\end{example}

\begin{definition}\label{defi2.4}
Let $p \in [2, \infty]$. A map
$$\rho : \FB := L^p(\FC_T) \rightarrow \RB$$
is called a {\it convex risk measure} if the following holds:
\begin{myenumerate}
\item (convexity) $\rho(\lambda F_1 + (1 - \lambda)F_2) \leq \lambda \rho(F_1) + (1 - \lambda) \rho (F_2)$ for all $\lambda \in [0,1]$ and all $F_1, F_2 \in \FB$.
\item (Monotonicity) If $F_1 \leq F_2$ then $\rho (F_1) \geq \rho (F_2)$.
\item (Translation invariance) $\rho (F + a) = \rho (F)-a$ for all $F\in \FB$ and all constants $a$.
\item (For convenience) $\rho(0) = 0$.
\end{myenumerate}
\end{definition}

\paragraph{Interpretation:} $\rho(F)$ is the amount that has to be added to the financial standing $F$ to make it ``acceptable''.
Note that, by (iii), $\rho(F + \rho (F)) = 0$.

\

We shall see that convex risk measures are related to BSDEs as follows: Let $g(t,z)$ be a concave function of $z$. For given $F \in \FB$ let $(Y^{(F)}_g(t) , Z^{(F)}_g(t))$ be the solution of the BSDE
\begin{equation}\label{eq2.11}
\begin{cases}
d Y^{(F)}_g(t) = - g(t, Z^{(F)}_g(t))dt + Z^{(F)}_g(t) dB(t) \; ; \; 0 \leq t \leq T \\
Y^{(F)}_g(T) = F.
\end{cases}
\end{equation}
Then
\begin{equation}\label{eq2.12}
\rho(F) := - Y^{(F)}_g(0)
\end{equation}
defines a convex risk measure.
The proof will be given in the next section, after we have studied BSDEs in a more general context with jumps.

\subsection{General BSDEs with jumps}\label{sec2.2}

Let $g(t,y,z,k,\omega) : [0,T] \times \RB \times \RB \times \RC \times \Omega \rightarrow \RB$ be a given function such that $(t,\omega) \rightarrow g(t,y,z,k,\omega)$ is $\FC_t$-predictable for each $y,z,k$.  Here $\RC$ is the set of functions $k \in L^2(\RB_0,\nu)$, where $\RB_0 := \RB \backslash \{0\}$. Let $F \in L^2(\FC_T)$. We seek a triple $(Y,Z,K) = (Y^{(F)}_g(t), Z^{(F)}_g(t), K^{(F)}_g(t))$ of stochastic processes such that $Y$ is a  c\`adl\`ag adapted process, $Z$ and $K$ are predictable and 
\begin{equation}\label{eq2.13}
\begin{cases}
dY(t) & = - g(t, Y(t), Z(t), K(t,\cdot),\omega) dt + Z(t) dB(t) \\ 
& \displaystyle + \int_\RB K(t,\zeta) \tN (dt, d\zeta) \; ; \; 0 \leq t \leq T \\
 Y(T) & = F.
 \end{cases}
 \end{equation}
 
 The process $g(t,y,z,k) = g(t,y,z,k,\omega)$ is called the {\it driver} of the BSDE \eqref{eq2.13}. We state the next result without proof. We refer to
  Tang and Li \cite{TL} (1994) and Quenez \& Sulem (2013), Theorem 2.3 for details.
 
 \begin{theorem}[Existence and uniqueness of solution of BSDE] \label{th2.5}

 Suppose the following holds:
 \begin{myenumerate}
 \item $ \displaystyle E \left[ \int_0^T g^2(t,0,0,0)dt \right] < \infty$
 \item $g$ is Lipschitz in $y,z,k$ a.s., i.e. there exists a constant $C > 0$ such that for all $y_i, z_i, k_i$ $|g(t,y_1, z_1, k_1) - g(t,y_2,z_2,k_2)| \leq C(|y_1-y_2|) + |z_1-z_2| + \| k_1-k_2\|)$ for a.e. $t,\omega$, where $\displaystyle \|k\|^2 = \int_\RB k^2(\zeta) \nu (d\zeta)$.
 \end{myenumerate}
 Then there exists a unique triple $(Y,Z,K)$  solution of \eqref{eq2.13} such that $Y$ is a c\`adl\`ag adapted  process with 
  $E(\sup_{0\leq t \leq T} |Y(t)|^2) <  \infty $  
  and $(Z,K)$ are predictable processes with
 \begin{equation}\label{eq2.14}
 E \left[ \int_0^T \left\{  Z^2(t) + \int_\RB K^2 (t,\zeta) \nu (d\zeta) \right\} dt \right] < \infty.
 \end{equation}
 \end{theorem}
 This result can be extended to the case when the terminal time  $T$  is  a stopping time $ \tau$ with values in $[0,T]$ and the terminal condition is  
 a random variable  $\xi$ in $ L^2({\cal F}_ \tau)$. 
In this case, $(Y^{(\xi, \tau)}, Z^{(\xi, \tau)}, K^{(\xi, \tau)})$ 
 is defined  as the unique solution of the BSDE with  driver $g(t,y, z,k) {\bf 1}_ { \{ t \leq  \tau  \}  }$ and terminal conditions ($T$, $\xi$). Note that 
$Y^{(\xi, \tau)}(t) = \xi, Z^{ (\xi, \tau)}(t)=K^{(\xi, \tau)}(t)=0$ for $t \geq  \tau$.
 \begin{lemma}\label{lem2.6}
 \begin{myenumerate}
 \item Suppose $(Y,Z,K)$ satisfies the BSDE \eqref{eq2.13}. Then
  \begin{equation}\label{eq2.15}
 Y(t) = E \left[ \left(\int_t^T g(s, Y(s), Z(s), K(s, \cdot)) ds + F \right) \mid \FC_t \right] \; ; \; 0 \leq t \leq T.
 \end{equation}
 \item Conversely, suppose the driver $g$ does not depend on $z$ and $k$ and that \eqref{eq2.15} holds. Then \eqref{eq2.13} holds.
 \end{myenumerate}
 \end{lemma}
 
 \dproof
 \paragraph{(i): \eqref{eq2.13} $\Rightarrow$ \eqref{eq2.15}}: Integrating \eqref{eq2.13} from $t$ to $T$ we get
$$
 Y(T) - Y(t)  = - \int_t^T g(s, Y(s), Z(s), K(s,\cdot))dt + \int_t^T Z(s) dB(s) + \int_t^T \int_\RB K(s,\zeta) \tN (ds, d\zeta).$$
 Taking conditional expectation and using that $\displaystyle t \rightarrow \int_0^t Z(s) dB(s)$ and  $\displaystyle t \rightarrow \int_t^T \int_\RB K(s,\zeta) \tN(ds, d\zeta)$ are martingales, we get \eqref{eq2.15}.
 \paragraph{(ii): \eqref{eq2.15} $\Rightarrow$ \eqref{eq2.13}}: Assume \eqref{eq2.15} holds and that $g(s,y,z,k) = g(s,y)$ does not depend on $z$ and $k$. Since $Y(T) = F$, we can write
 \begin{align}\label{eq2.16}
 Y(t) & = E \left[ \int_0^T \{g(s, Y(s))dt + F - \int_0^t g(s, Y(s))ds \} \mid \FC_t \right] \nonumber \\
 & = M(t) - \int_0^t g(s, Y(s))ds,
 \end{align}
 where $M(t)$ is the $L^2$-martingale
 $$M(t) = E \left[ \int_0^T \{g(s, Y(s))ds + F\} \mid \FC_t \right].$$
 By the martingale representation theorem  for It\^o-L\'evy process (see e.g. \cite{L}) there exists $Z(t)$ and $K(t, \zeta)$ such that
 \begin{equation}\label{eq2.17}
 M(t) = M(0) + \int_0^t Z(s) dB(s) + \int_0^t \int_\RB K(s, \zeta) \tN (ds,d\zeta) \; ; \; 0 \leq t \leq T.
 \end{equation}
 Substituting \eqref{eq2.17} into \eqref{eq2.16} and taking differentials, we get \eqref{eq2.13}.
 \fproof
 \subsection{Linear BSDEs}
 
 There is no solution formula for the general BSDE \eqref{eq2.13}. However, in the {\it linear} case we get the following:
 
 \begin{theorem}\label{th2.7}
 Let $\alpha, \beta, \gamma$ be bounded predictable processes, $F \in L^2(\FC_T)$ and $\varphi$ predictable with $\displaystyle E \left[ \int_0^T \varphi^2(t)dt\right] < \infty$. Assume  $\gamma > -1$ a.s. Then the unique solution $(Y,Z,K)$ of the linear BSDE
 \begin{equation}
 \label{eq2.18}
 \begin{cases}
 dY(t)  \displaystyle = - \left[ \varphi(t) + \alpha(t) Y(t) + \beta(t) Z(t) + \int_\RB \gamma(t,\zeta) K(t,\zeta) \nu (d \zeta)\right]dt \\
  \quad \displaystyle + Z(t) dB(t) + \int_\RB K(t,\zeta) \tN (dt, d\zeta) \; ; \; 0 \leq t \leq T \\
  Y(T) = F
  \end{cases}
  \end{equation}
  is given by
  \begin{equation}\label{eq2.19}
  Y(t) = E \left[ \{\frac{\Gamma(T)}{\Gamma(t)} F + \int_t^T \frac{\Gamma(s)}{\Gamma(t)} \varphi(s) ds\} \mid \FC_t \right] \; ; \; 0 \leq t \leq T
  \end{equation}
  where
  \begin{equation}\label{eq2.20}
  \begin{cases}
 \displaystyle  d\Gamma(t) = \Gamma(t^-) \left[ \alpha(t) dt + \beta(t) dB(t) + \int_\RB \gamma(t,\zeta) \tN (dt,d\zeta)\right] \; ; \; t \geq 0 \\
 \Gamma(0) = 1
 \end{cases}
 \end{equation}
 i.e.
 \begin{align}\label{eq2.21}
 \Gamma(t) & = \exp \left( \int_0^t \beta(s) dB(s) + \int_0^t \left\{ \alpha(s) - \frac{1}{2} \beta^2(s) \right\} ds \right.\nonumber \\
  & + \int_0^t \int_\RB ln (1 + \gamma (s,\zeta))\tN (ds,d\zeta) + \int_0^t \int_\RB \{ ln (1 + \gamma(s,\zeta)) - \gamma(s,\zeta)\} \nu (d\zeta) ds.
  \end{align}
 \end{theorem}
 
 \dproof (Sketch). By the It\^o formula,
 
 \begin{align*}
 d&(\Gamma(t)Y(t)) = \Gamma(t^-) dY(t) + Y(t^-)d \Gamma(t) +d [ \Gamma Y](t) \\
 & = \Gamma(t^-) \left[ -\left(\varphi(t) + \alpha(t) Y(t) + \beta(t) Z(t) + \int_\RB \gamma(t,\zeta) K(t,\zeta) \nu(d\zeta)\right)\right. dt \\
 & \left. \quad + Z(t) dB(t) + \int_\RB K(t,\zeta) \tN (dt, d\zeta)\right]\\
 & \quad + Y(t) \Gamma(t^-) \left[ \alpha(t) dt + \beta(t) dB(t) + \int_\RB \gamma(t,\zeta) \tN (dt,d\zeta)\right] \\
 & \quad + Z(t) \Gamma(t) \beta(t) dt + \int_\RB K(t,\zeta) \Gamma(t^-) \gamma(t,\zeta) N(dt,d\zeta) \\
 & = - \Gamma(t) \varphi(t)dt + \Gamma(t) (Z(t) + \beta(t) Y(t)) dB(t) \\
  & \quad + \int \Gamma(t^-) K(t,\zeta) (1 + \gamma(t,\zeta)) \tN (dt,d\zeta).
 \end{align*}
Hence $\displaystyle \Gamma(t)Y(t) + \int_0^t \Gamma(s) \varphi(s)ds$ is a martingale and therefore
$$\Gamma(t) Y(t) + \int_0^t \Gamma(s) \varphi(s) ds = E \left[ \{\Gamma(T) Y(T) + \int_0^T \Gamma(s) \varphi(s)ds\} \mid \FC_t \right]$$
i.e.
$$
\Gamma(t) Y(t) = E \left[ \{\Gamma(T) F + \int_t^T \Gamma(s) \varphi(s) ds\} \mid \FC_t\right],$$
as claimed.
 \fproof
 
 \begin{example}\label{ex2.8} Let us apply Theorem \ref{th2.7} to solve the BSDE \eqref{eq2.6}-\eqref{eq2.3}:
 In this case
 $$d \Gamma(t) = \Gamma(t) \left[ r(t) dt + \frac{\mu(t) - r(t)}{\sigma(t)} dB(t)\right]  \; ; \; \Gamma(0) = 1$$
 i.e.
 $$\Gamma(t) = \exp \left( \int_0^t \frac{\mu(s) - r(s)}{\sigma(s)} dB(s) + \int_0^t \left\{ r(s) - \frac{1}{2} \left( \frac{\mu(s) - r(s)}{\sigma(s)}\right)^2 \right\}ds \right)$$
 and we get
 $$Y(t) = \frac{1}{\Gamma(t)} E [ F \Gamma(T)| \FC_t].$$
 Using Malliavin calculus we can write
 $$Z(t) =D_{t^-} Y(t) \left( := \lim_{s \rightarrow t^-} D_s Y(t)\right)$$
 and this gives the replacing portfolio
 $$\pi(t) =\frac{Z(t)}{Y(t)\sigma(t)} = \frac{D_tY(t)}{Y(t) \sigma(t)}.$$
 Here $D_t$ denotes the Malliavin derivative at $t$ (with respect to Brownian motion). See e.g. Di Nunno et al. (2009).
 
 \end{example}
 
 \subsection{Comparison theorems}\label{sec2.4}

\begin{lemma}\label{lem2.9}
Let $\alpha, \beta, \gamma, F$ be as in Theorem \ref{th2.7}. Suppose $(Y(t), Z(t), K(t,\cdot))$ satisfies the linear backward stochastic inequality
\begin{equation}\label{eq2.22} 
\begin{cases}
\displaystyle dY(t) = - h(t) dt + Z(t) dB(t) + \int_\RB K(t,\zeta) \tN (dt, d\zeta) \; ; \; 0 \leq t \leq T \\
Y(T) \geq F
\end{cases}
\end{equation}
where $h(t)$ is a given $\FC_t$-adapted process such that
\begin{equation}\label{eq2.23}
h(t) \geq  \alpha(t) Y(t) + \beta(t) Z(t) + \int_\RB \gamma(t,\zeta) K(t,\zeta)\nu (d \zeta).
\end{equation}
Then
\begin{equation}\label{eq2.24}
Y(t) \geq E \left[ \frac{\Gamma(T)F}{\Gamma(t)} \mid \FC_t \right] \; ; \; 0 \leq t \leq T
\end{equation}
where $\Gamma(t)$ is given by \eqref{eq2.20}-\eqref{eq2.21}.
\end{lemma}

\dproof
By the It\^o formula we have
\begin{align*}
d(\Gamma(t) Y(t))& = \Gamma(t^-) \left[ - h(t) dt + Z(t) dB(t) + \int_\RB K(t, \zeta) \tN (dt, d\zeta) \right] \\
 &\quad + Y(t^-) \Gamma(t^-) \left[ \alpha(t) dt + \beta(t) dB(t) + \int_\RB \gamma(t,\zeta) \tN (dt, d\zeta) \right] \\
 &\quad + \Gamma(t) \beta(t) Z(t) dt + \int_\RB \Gamma(t^-) \gamma(t,\zeta)K(t,\zeta) N (dt, d\zeta) \\
 & \leq \Gamma(t) \left[  - \alpha(t) Y(t) - \beta(t) Z(t) - \int_\RB \gamma(t,\zeta) K(t,\zeta) \nu (d \zeta) \right]dt \\
 & \quad + \Gamma(t) Z(t) dB(t) + \Gamma(t^-) \int_\RB K(t,\zeta) \tN (dt, d \zeta) \\
 & \quad + Y(t^-)  \Gamma (t^-) \left[ \alpha(t) dt + \beta(t) dB(t) + \int_\RB \gamma(t,\zeta) \tN (dt, d\zeta) \right] \\
 & \quad + \Gamma(t) \beta(t) Z(t) dt + \int_\RB \Gamma(t^-) \gamma(t,\zeta) K(t,\zeta) N (dt,d\zeta) \\
 & = d M(t),
\end{align*}
where 
$$M(t) := \int_0^t \Gamma (s) Z(s) dB(s) + \int_0^t \int_\RB \Gamma(s^-) \gamma(s,\zeta) K(s,\zeta) \tN (ds, d\zeta)$$ 
is a martingale.
Hence
$$\Gamma(T) Y(T) - \Gamma(t) Y(t) \leq  M(T) - M(t).$$
Taking conditional expectation this gives
$$\Gamma(t) Y(t) \geq E \left[ \Gamma(T) F  \mid \FC_t \right].$$
\fproof
\begin{corollary}\label{cor2.10}
Let $Y, Z, K$ be as in Lemma \ref{lem2.9}. Suppose that
$$ F \geq 0 \text{ a.s. }$$
Then
$$Y(t) \geq 0 \text{ for a.a. } t, \omega.$$
\end{corollary}
\dproof
Apply Lemma \ref{lem2.9}.
\fproof

In the following we assume that
$$g_i(t,y,z,k,\omega) : [0,T] \times \RB \times \RB \times \RC \times \Omega \rightarrow \RB \; ; \; i = ,2$$
 are given $\FC_t$-predictable processes  satisfying (i)-(ii)  in Theorem \ref{th2.5}. We assume that $g_2(t,y,z,k,\omega)$ is Lipschitz continuous with respect to $y,z,k$, uniformly in $t,\omega$.
We also assume that there exists a bounded predictable process $\theta(t,\zeta)$ independent of $y$ and $z$ 
such that $dt\otimes dP \otimes \nu(du)$-a.s.\,,
\begin{equation}\label{robis}
\theta(t, \zeta) \geq -1  \;\; \text{ and }
\;\; |\theta(t, \zeta) | \leq \psi (\zeta),
\end{equation}
where $\psi$ $\in$ $L^2_{\nu}$, and such that 
\begin{equation}\label{eq2.25}
g_2(t,y,z,k_1(\cdot)) - g_2(t,y,z,k_2(\cdot)) \geq \int_\RB \theta(t,\zeta) (k_1(\zeta) - k_2(\zeta)) \nu (d\zeta)
\end{equation}
for all $t,y,z$.

We are now ready to state and prove a comparison theorem for BSDEs with jumps. For a stronger version see \cite{QS}.

\begin{theorem}(Comparison theorem for BSDEs with jumps)\label{th2.10}
Suppose we have 2 process triples $(Y_1,Z_1,K_1)$ and $(Y_2, Z_2, K_2)$, such that
\begin{equation}\label{eq2.26}
\begin{cases}
dY_i(t) & = - g_i(t,Y_i(t), Z_i(t), K_i(t,\cdot)) dt + Z_i(t) dB(t) \\
& \displaystyle + \int_\RB  K_i(t,\zeta) \tN (dt, d\zeta) \; ; \; 0 \leq t \leq T \\
Y_i(T) & = F_i
\end{cases}
\end{equation}
for $i = 1,2$.
 where  $F_i  \in L^2(\FC_T).$
Assume that
\begin{equation}
\label{eq2.27}
g_1(t,Y_1(t), Z_1(t), K_1(t, \cdot)) \leq g_2(t, Y_1(t), Z_1(t), K_1(t, \cdot)) \; ; \; t \in [0,T]
\end{equation}
and
\begin{equation}\label{eq2.28}
F_1 \leq F_2 \text{ a.s.}
\end{equation}
Then
\begin{equation}\label{eq2.29}
Y_1(t) \leq Y_2(t) \text{ for a.a. } (t,\omega) \in [0,T] \times \Omega.
\end{equation}
\end{theorem}

\dproof
Put
$$\tY(t) = Y_2(t) - Y_1(t), \tZ(t) = Z_2(t) - Z_1(t), \tK(t,\zeta) = K_2(t,\zeta) - K_1(t,\zeta).$$
Then
\begin{align*}
d\tY(t) & = - [g_2(t,Y_2(t),Z_2(t),K_2(t,\cdot)) - g_1(t, Y_1(t), Z_1(t), K_1(t, \cdot))] dt \\
 & + \tZ(t) dB(t) + \int_\RB \tK(t,\zeta) \tN(dt, d\zeta) \; ; \; 0 \leq t \leq T.
 \end{align*}
 Note that
 \begin{align*}
 g_2&(t,Y_2(t), Z_2(t), K_2(t,\cdot)) - g_1(t, Y_1(t), Z_1(t), K_1(t, \cdot)) \\
 & = g_2(t, Y_2(t), Z_2(t), K_2(t, \cdot)) - g_2(t, Y_1(t), Z_2(t), K_2(t, \cdot)) \\
 & \quad + g_2(t, Y_1(t), Z_2(t), K_2(t,\cdot)) - g_2(t, Y_1(t), Z_1(t), K_2(t, \cdot)) \\
 & \quad + g_2 (t,Y_1(t), Z_1(t), K_2(t, \cdot)) - g_2(t, Y_1(t), Z_1(t), K_1(t, \cdot)) \\
 & \quad + g_2(t, Y_1(t), Z_1(t), K_1(t, \cdot)) - g_1(t, Y_1(t), Z_1(t), K_1(t, \cdot)) \\
 & = \varphi(t) + \alpha(t) \tY(t) + \beta(t) \tZ(t) + \int_\RB \theta(t,\zeta) \tK(t,\zeta) \nu(d \zeta),
 \end{align*}
 by \eqref{eq2.25}, where
 $$\varphi(t) := g_2(t,Y_1(t), Z_1(t), K_1(t, \cdot)) - g_1(t, Y_1(t), Z_1(t), K_1(t, \cdot)) \geq 0,$$
 $$\alpha(t) := \frac{g_2(t,Y_2(t), Z_2(t), K_2(t,\cdot)) - g_2(t, Y_1(t), Z_2(t), K_2(t, \cdot))}{\tY(t)} \chi_{\tY(t) \neq 0} \tY(t)$$
 and
 $$\beta(t) := \frac{g_2(t, Y_1(t), Z_2(t), K_2(t, \cdot)) - g_2(t, Y_1(t), Z_1(t), K_2(t, \cdot))}{\tZ(t)} \chi_{\tZ(t) \neq 0} \tZ(t).$$
 Combining the above we get
 $$\begin{cases}
 d \tY(t) & = \displaystyle - h(t) dt + \tZ (t) dB(t) + \int_\RB \tK(t,\zeta) \tN (dt, d\zeta) \; ; \; 0 \leq t \leq T \\
 \tY(T) & = F_2 - F_1 \geq 0
 \end{cases}
 $$
 where
 $$h(t) \geq  \alpha(t) \tY(t) + \beta (t) \tZ(t) + \int_\RB \theta(t,\zeta) \tK(t,\zeta) \nu (d\zeta).$$
 By Corollary \ref{cor2.10} it follows that $\tY(t) \geq 0$ for all $t$, i.e. $Y_1(t) \leq Y_2(t)$ for all $t$.
\fproof
\subsection{Convex risk measures, recursive utilities and BSDEs}\label{sec2.5}
We now have the machinery we need to verify the connection between risk measures, recursive utilities and BSDEs mentioned in Section \ref{sec2.1}.
Motivated by Lemma \ref{lem2.6} we now extend the definition of recursive utility given in Example \ref{ex2.2} to the following:
We call a process $c(t)$ a {\it consumption} process if $c(t)$ is predictable and $c(t) \geq 0$ for all $t$, a.s. The set of all consumption processes is denoted by $\CC$.
\begin{definition}\label{defi2.12}
  Let $g(t,y,z,k,c) : [0,T] \times \RB \times \RB \times \RC \times \CC \rightarrow \RB$ be a process satisfying the conditions of Theorem \ref{th2.5} for each given $c \in \CC$. Suppose
\begin{equation}\label{eq2.30} (y,z,k,c) \rightarrow g(t,y,z,k,c) \text{ is concave for all } t.
\end{equation}
Let $(Y^{(F)}_g, Z_g^{(F)}, K_g^{(F)})$ be the unique solution of the BSDE \eqref{eq2.13}. Then we define
\begin{equation}\label{eq2.31}
U_g(c) = Y_g^{(F)}(0)
\end{equation}
to be the recursive utility of $c$ with terminal payoff $F$.
\end{definition}

\begin{theorem}\label{th2.13} Suppose $g(t,z,k) : [0,T] \times \RB \times \RB \times \RC \rightarrow \RB$ satisfies the conditions in Definition \ref{defi2.12}, but now $g$ does not depend on $y$ or $c$.  Assume $g$ satisfies hypothesis \eqref{eq2.25}. Define
\begin{equation}\label{eq2.32}
\rho_g(F) = - Y_g^{(F)}(0).
\end{equation}
Then $\rho_g$ is a convex risk measure.
\end{theorem}

\dproof

 We must verify that $\rho_g$ satisfies the properties (i)-(iii) in Definition \ref{defi2.4}:
\begin{myenumerate}
\item (Convexity). Fix $\lambda \in (0,1)$ and let $F,G \in L^2(\FC_T)$. We want to prove that
$$\rho_g (\lambda F + (1- \lambda)G) \leq \lambda \rho_g(F) + (1 - \lambda) \rho_g(G)$$
i.e.
$$- Y^{(\lambda F + (1- \lambda)G)} (0) \leq \lambda (- Y^{(F)}(0)) + (1 - \lambda) (- Y^{(G)}(0)).$$
Let $(\hY, \hZ, \hK)$ be the solution of the BSDE
$$\begin{cases}
d\hY(t) & = \displaystyle - g(t, \hZ(t), \hK(t, \cdot)) dt + \hZ(t) dB(t) + \int_\RB \hK(t,\zeta) \tN(dt, d\zeta) \; ; \; 0 \leq t \leq T \\
\hY(T) & = \lambda F + (1 - \lambda) G
\end{cases}
$$
and put
\begin{align*}\tY(t)& = \lambda Y^{(F)}(t) + (1 - \lambda) Y^{(G)}(t), \\
\tZ(t)& = \lambda Z^{(F)}(t) + (1 - \lambda) Z^{(G)}(t), \\
\tK(t,\zeta) &= \lambda K^{(F)}(t,\zeta) + (1 - \lambda) K^{(G)}(t,\zeta).\end{align*}
Then
$$\begin{cases}
d \tY(t) & = - [h(t) + g(t, \tZ(t), \tK(t, \cdot))]dt \\
& \displaystyle + \tZ (t) dB(t) + \int_\RB \tK(t,\zeta) \tN (dt, d\zeta) \; ; \; 0 \leq t \leq T \\
\tY(T) & = \lambda F + (1 - \lambda) G,
\end{cases}
$$
where
\begin{align*}
h(t) & = \lambda g(t, Z^{(F)}(t), K^{(F)}(t)) + (1 - \lambda) g(t, Z^{(G)}(t), K^{(G)}(t,\cdot)) \\
& - g(t, \tZ(t), \tK(t,\cdot)) \leq 0 \text{ since $g$ is concave.}
\end{align*}
By the comparison theorem (Theorem \ref{th2.10}) we conclude that
$$\tY(t) \leq \hY(t) \; ; \; 0 \leq t \leq T.$$
In particular, choosing $t = 0$ we get
\begin{align*}
\rho_g(\lambda F + (1-\lambda)G) & = - \hY(0) \leq - \tY(0) = - \lambda Y^{(F)}(0) - (1 - \lambda) Y^{(G)}(0) \\
& = \lambda \rho_g (F) + (1 - \lambda) \rho_g (G).
\end{align*}
\item (Monotonicity) If $F_1 \leq F_2$, then $Y^{(F_1)}(t) \leq Y^{(F_2)}(t)$ by the comparison theorem. Hence
$$\rho_g(F_2) = - Y^{(F_2)}(0) \leq - Y^{(F_1)}(0) = \rho_g (F_1),$$
as required.
\item (Translation invariance) If $F \in L^2(\FC_T,P)$ and $a \in \RB$ is constant, then we check easily that $Y^{(F+a)}(t) = Y^{(F)}(t) + a$. Hence
$$\rho(F+a) = - Y^{(F+a)}(0) = - Y^{(F)}(0) - a = \rho_g (F) -a.$$
\end{myenumerate}
\fproof
\paragraph{Dynamic risk measures.}

We now discuss an extension of the (static) risk measure $\rho$ in Definition \ref{defi2.4} to a {\em dynamic} 
risk measure $\rho_t;  \;\; 0 \leq t \leq T.$

\begin{definition}\label{defi2.14} \rm
A {\em dynamic risk measure} is a map $\rho$ which to each bounded stopping time $\tau$ and each $\xi \in L^2({\cal F}_\tau)$ 
assigns an adapted c\`adl\`ag process $(\rho_t(\xi, \tau))_{\{0 \leq t \leq \tau\}}$  which is non-increasing, translation invariant and 
{\em consistent}, in the sense that
\begin{equation}\label{con}
 \forall t \leq S, \rho_t(\xi, \tau) = \rho_t (-\rho_S(\xi, \tau), S) \text{ a.s.} 
 \end{equation}
for all stopping times $S \leq \tau$. 

Moreover we say that the risk measure satisfies 
\begin{itemize}
\item the {\em zero-one law} property if \\
$\rho_t({\bf 1}_A \xi, T) = {\bf 1}_A \rho_t (\xi, T)$ $a.s$ for $t \leq T$, $A \in {\cal F}_t$, and $\xi$ $\in$ $L^2({\cal F}_T)$.
\item the  {\em no arbitrage} property if \\
\big\{$\xi^1 \geq \xi^2$   a.s. and $ \rho_t(\xi^1,\tau) = \rho_t(\xi^2, \tau)$ a.s. on some $A \in  {\cal F}_t$,   $t \leq \tau$\big\}
$\Longrightarrow$ \big\{$ \xi^1 = \xi^2$ a.s. on $A$\big\} .
\end{itemize}
\end{definition}
%
%
%
%
A  natural way to construct dynamic risk measures is by means of BSDEs as follows:\\
Let $g$ be a Lipschitz driver, which does not depend on $y$ and  such that $E \left[ \int_0^T g^2(t,0,0)dt \right] < \infty$. 
We assume that $g$ satisfies \eqref{eq2.25} -\eqref{robis} with $\theta(t, \zeta) > -1$.
 For a given stopping time $\tau \leq T$ and  $\xi \in L^2({\cal F}_T)$, define
 the  functional:
  \begin{equation}\label{definition}
  \rho^g_t(\xi,\tau)  :=  -Y_g^{(\xi)}(t), \,\, \,\,\,0\leq t \leq \tau,
  \end{equation}
 where $Y_g^{(\xi)}$ denotes the solution  of the BSDE with terminal condition $\xi$ and terminal time  $\tau$.
 Then  $\rho^g$ defines  a {\em dynamic risk measure} in the sense of Definition \ref{defi2.14}.
To see this, we note that the consistency \eqref{con} follows  from the flow property of BSDEs (see \cite{QS}).

Moreover, the no-arbitrage property follows from the strict comparison theorem for BSDEs.
We also note that if  $g(t,0,0) = 0$, then the zero-one law holds.
The dynamic risk measure is convex if $g$ is concave. 

It is natural to ask about the converse: 
When can  a dynamic risk-measure be represented by a BSDE with jumps? The following proposition gives an answer. 
\begin{theorem} 
Let $\rho$ be  a dynamic risk measure 
 satisfying the zero-one law and the no arbitrage property. Moreover, suppose that $\rho$ satisfies the  {\em ${\cal E}^{C,C^1}$-domination} property: \\
there exists  $C >0$ and $-1 < C_1 \leq 0$ such that 
\begin{equation}\label{risk}
\rho_t(\xi +\xi ^{\prime}, T) - \rho_t(\xi , T) \geq - Y_t ^{C, C_1} (\xi ^{\prime},T),
\end{equation}
for any $\xi,\xi ^{\prime}$ $\in$ $L^2({\cal F}_T)$, where $Y_t ^{C, C_1} (\xi ^{\prime},T)$ is the solution of the BSDE associated with terminal condition $\xi^{\prime}$ and driver 
$f_{C, C_1}(t, \pi, \ell) := C |\pi|  + C \int_{\R^*} (1 \wedge |u|) \ell^+(u) \nu(du) - C_1 
\int_{\R^*} (1 \wedge |u|) \ell^-(u) \nu(du)$.
Then,  there exists a Lipschitz driver $g(t, \pi, \ell)$ such that $\rho=\rho^g$, that is,  $\rho$ is the dynamic risk measure induced by a BSDE with jumps with  driver $g(t, \pi, \ell)$ .  
\end{theorem}

 For the proof,  we refer to \cite{R}. 
Additional properties of dynamic risk measures induced by BSDEs and dual representation in the convex case can be found in \cite{QS}.

\section{Maximum principles for optimal control of coupled systems of FBSDEs}\label{sec3}

In view of Definition \ref{defi2.12} and Theorem \ref{th2.13}, we see that recursive utility maximization or risk minimization problems lead to problems of optimal control of coupled systems of forward-backward stochastic differential equations (FBSDEs). In this section we study such control problems. For simplicity we only handle the 1-dimensional case.

Consider the following stochastic control problem for a system of coupled forward-backward stochastic differential equations (FBSDEs):

\noindent (Forward system)
\begin{equation}\label{eq3.1}
\begin{cases}
dX(t) & = b(t, X(t), Y(t), Z(t), K(t,\cdot), u(t), \omega) dt \\
 & + \sigma(t, X(t), Y(t), Z(t), K(t,\cdot), u(t), \omega) dB(t) \\
  & \displaystyle + \int_\RB \gamma(t, X(t) Y(t), Z(t) K(t,\cdot),u(t), \omega, \zeta) \tN (dt, d\zeta) \; ;\; t \geq 0 \\
  X(0) & = x \in \RB.
\end{cases}
\end{equation}

\noindent (Backward system)
\begin{equation}\label{eq3.2}
\begin{cases}
dY(t) = & - g(t, X(t), Y(t), Z(t), K(t,\cdot), u(t), \omega)dt \\
& \displaystyle + Z(t) dB(t) + \int_\RB K(t,\zeta) \tN (dt, d\zeta) \; ; \; 0 \leq t \leq T \\
Y(T) & = h(X(T)).
\end{cases}
\end{equation}
Here $T > 0$ is fixed (finite) constant. 
Let
$\GB:=\{\GC_t\}_{0 \leq t \leq T}$ be a given \emph{subfiltration} of $\FB:=\{\FC_t\}_{0 \leq t \leq T},$  i.e. $\GC_t \subseteq \FC_t$ for all $t$. We assume that also $\GB$ satisfies the usual conditions. We can interpret $\GC_t$ as the information available to the controller at time $t$.\\
Let $\UB$ be a given open convex subset of $\RB$ and let $\AC_\GB$ be a given family of admissible controls, consisting of all $\GB$-predictable processes $u=u(t)$ with values in $\UB$.

The {\it performance functional} is given by
\begin{align}\label{eq3.3}
J(u) & = E \left[ \int_0^T f(t,X(t), Y(t), Z(t), K(t,\cdot), u(t), \omega) dt 
  + \varphi (X(T), \omega)\right] \nonumber \\
  & + \psi (Y(0)) \; ; \; u \in \AC_\GB,
\end{align}
We want to find $u^* \in \AC_\GB$ such that
\begin{equation}\label{eq3.4}
\sup_{u \in \AC_\GB} J(u) = J(u^*).
\end{equation}
We make the following assumptions:
\begin{align}\label{eq3.5}
f \in C^1 \text{ and } E \left[ \int_0^T | \nabla f |^2 (t)dt \right] < \infty, \\
b,\sigma, \gamma \in \CC^1 \text{ and } E \left[ \int_0^T (| \nabla b|^2 + |\nabla \sigma|^2 + \| \nabla \gamma\|^2)(t)dt \right] < \infty, \label{eq3.6}
\end{align}
where $\displaystyle \|\nabla \gamma(t,\cdot)\|^2 = \int_\RB \gamma^2(t,\zeta) \nu (d \zeta)$,
\begin{equation}\label{eq3.7}
g \in \CC^1 \text{ and } E \left[ \int_0^T | \nabla g |^2 (t) dt \right] < \infty,
\end{equation}
\begin{equation}\label{eq3.8}
h, \varphi, \psi \in \CC^1 \text{ and } E[ \varphi'(X(T))^2 + h'(X(T))^2] < \infty\end{equation}
for all $u \in \AC_\GB$.
Let $\RC$ denote the set of all functions $k : \RB_0 \rightarrow \RB$ where $\RB_0 = \RB \backslash \{0\}$.

The \emph{Hamiltonian}
$$H:[0,T] \times \RB \times \RB \times \RB \times \RC \times \UB \times \RB \times \RB \times \RB \times \RC \times \Omega \rightarrow \RB$$
associated to the problem \eqref{eq3.4} is defined by
\begin{align}\label{eq3.9}
H(t,x,y,z,k,u,\lambda,p,q,r,\omega)& = f(t,x,y,z,k,u,\omega) + g(t,x,y,z,k,u,\omega) \lambda + b(t,x,y,z,k,u,\omega)p \nonumber \\
& \quad + \sigma(t,x,y,z,k,u,\omega)q + \int_\RB \gamma(t,x,y,z,k,u,\zeta,\omega) r(t,\zeta)\nu (d\zeta).
\end{align}
Here $\lambda, p,q,r$ represent adjoint variables (see below).

For simplicity of notation the dependence on $\omega$ is suppressed in the following.

We assume that $H$ is Fr\'echet differentiable $(C^1)$ in the variables $x,y,z,k,u$ and that the Fr\'echet derivative $\nabla_k H$ of $H$ with respect to $k \in \RC$ as a random measure is absolutely continuous with respect to $\nu$, with Radon-Nikodym derivative $\displaystyle \frac{d \nabla_k H}{d \nu}$. Thus, if $\langle \nabla_k H, h \rangle$ denotes the action of the linear operator $\nabla_k H$ on the function $h \in \RC$ we have
\begin{equation}\label{3.10}
\langle \nabla_k H, h \rangle = \int_\RB h(\zeta) d \nabla_k H(\zeta) = \int_\RB h(\zeta) \frac{d \nabla_kH(\zeta)}{d \nu (\zeta)} d \nu (\zeta).
\end{equation}
We let $m$ denote Lebesgue measure on $[0,T]$.
For $u \in \AC_\GB$ we let $(X^u(t), Y^u(t), Z^u(t), K^u(t,\cdot))$ be the associated solution of the coupled system \eqref{eq3.1}-\eqref{eq3.2}. We assume that for $u \in \AC_\GB$ these solutions exist and are unique and satisfy
\begin{equation}\label{eq3.10}
E \left[ \int_0^T \left\{ |X^u(t)|^2 + |Y^u(t)|^2 + Z^u(t)|^2 + \int_\RB | K^u(t,\zeta)|^2 \nu (d \zeta)\right\} dt \right] < \infty.
\end{equation}

The associated FB system for the adjoint processes $\lambda(t), (p(t), q(t), r(t,\cdot))$ is
\begin{equation}\label{eq3.11}
\begin{cases}
d\lambda(t)& = \displaystyle \frac{\partial H}{\partial y}(t) dt + \frac{\partial H}{\partial z}(t) dB(t) + \int_\RB \frac{d \nabla_k H}{d \nu}(t,\zeta) \tN(dt, d\zeta) \; ; \; 0 \leq t \leq T \\
\lambda(0) & = \psi'(Y(0))
\end{cases}
\end{equation}

\begin{equation}\label{eq3.12}
\begin{cases}
dp(t) & = \displaystyle - \frac{\partial H}{\partial x}(t)dt + q(t)dB(t) + \int_\RB r(t,\zeta) \tN (dt, d\zeta) \; ; \; 0 \leq t \leq T \\
p(T) & = \varphi'(X(T)) + \lambda(T) h'(X(T)).
\end{cases}
\end{equation}
Here and in the following we are using the abbreviated notation
$$\frac{\partial H}{\partial y}(t) = \left[\frac{\partial}{\partial y} H(t, X(t), y, Z(t), K(t,\cdot),u(t))\right]_{y = Y(t)}  \text{ etc.}$$

We first formulate a sufficient maximum principle. It is stronger than the corresponding result in e.g. \O ksendal \& Sulem (2012) because of our  weaker growth conditions here.

\begin{theorem}(Strengthened sufficient maximum principle)\label{th3.1}

Let $\hu \in \AC_\GB$ with corresponding solutions $\hX(t), \hY(t), \hZ(t), \hK(t,\cdot), \hla(t), \hp(t), \hq(t),\hr(t,\cdot)$ of equations \eqref{eq3.1}-\eqref{eq3.2}, \eqref{eq3.11} and \eqref{eq3.12}. Assume the following:
\begin{equation}\label{eq3.13}
\text{The functions } x \rightarrow h(x), x \rightarrow \varphi(x) \text{ and } x \rightarrow \psi(x) \text{ are concave}
\end{equation}
\begin{align}\label{eq3.14}
&\text{(The Arrow condition). The function } \nonumber\\
&\HC(x,y,z,k) := \ess_{v \in \UB} E [H(t,x,y,z,k,v, \hla(t), \hp(t), \hq(t), \hr(t,\cdot)) \mid \GC_t] \nonumber\\
&\text{is concave for all $t$, a.s.}
\end{align}
\begin{align}\label{eq3.15}
&\text{(The conditional maximum principle)} \nonumber\\
& \ess_{v \in \UB} E[H(t,\hX(t), \hY(t), \hZ(t), \hK(t,\cdot),v,\hla(t), \hp(t), \hq(t), \hr(t,\cdot)) \mid \GC_t] \nonumber\\
& \quad = E[H(t,\hX(t), \hY(t), \hZ(t), \hK(t,\cdot), \hu(t), \hla(t), \hp(t), \hq(t), \hr(t,\cdot)) \mid \GC_t] \; ; \; t \in [0,T]
\end{align}
\begin{equation}\label{eq3.16}
\left\| \frac{d \nabla_k \hH(t,.)}{d \nu} \right\| < \infty \text{ for all } t \in [0,T].
\end{equation}
Then $\hu$ is an optimal control problem the problem \eqref{eq3.4}.
\end{theorem}
\dproof
Define a sequence of stopping times $\tau_n \; ; \; n = 1,2, \ldots,$ as follows
\begin{align}\label{eq3.17}
\tau_n & = \inf\{ t > 0 \; ; \; \max \{ | \hp(t)|, | \sigma(t) - \hsi(t)|,\| \gamma(t,\cdot) - \hga(t,\cdot)\|, |X(t) - \hX(t)|, |\hq(t)|, \nonumber \\
  & \|\hr(t,\cdot)\|, |Y(t) - \hY(t)|, \left| \frac{\partial \hH}{\partial z}(t)\right|, \left\| \frac{d\nabla_k \hH}{d\nu}(t,\cdot)\right\|, | \hla(t)|, |Z(t) - \hZ(t)|\nonumber \\
   & \|K(t,\cdot) - \hK(t,\cdot)\| \} \geq n \} \wedge T.
\end{align}
Then note that $\tau_n \rightarrow T$ as $n \rightarrow \infty$ and
\begin{align}\label{eq3.18}
E &\left[ \int_0^{\tau_n} \hp(t) \left\{  (\sigma(t) - \hsi(t)) dB(t) + \int_\RB (\gamma(t,\zeta) - \hga(t,\zeta))\tN (dt d\zeta)\right\}\right] \nonumber \\
& = E \left[ \int_0^{\tau_n} (X(t) - \hX(t)) \left\{ \hq(t) dB(t) + \int_\RB \hr(t,\zeta) \tN  (dt,d\zeta)\right\}\right] \nonumber \\
& = E \left[ \int_0^{\tau_n}(Y(t^-) - \hY(t^-)) \left\{ \frac{\partial \hH}{\partial z}(t) dB(t) + \int_\RB \frac{ d\nabla_k \hH}{d \nu}(t,\cdot) \tN(dt, d\zeta) \right\}\right] \nonumber \\
&= E \left[ \int_0^{\tau_n} \hla(t) \left\{ (Z(t) - \hZ(t))dB(t) + \int_\RB (K(t,\zeta) - \hK(t,\zeta))\tN(dt,d\zeta) \right\}\right] \nonumber \\
& = 0 \text{ for all $n$.}
\end{align}
Except for the introduction of these stopping times, the rest of the proof follows the proof in \O ksendal \& Sulem (2012). For completeness we give the details: \\
Choose $u \in \AC_\GB$ and consider
$$J(u) - J(\hu) = J_1 + J_2 + J_3,$$
where
$$J_1 = E \left[ \int_0^T \{ f(t) - \hf(t)\} dt \right],\; J_2 = E [ \varphi(X(T)) - \varphi(\hX(T))], \; J_3 = \psi (Y(0)) - \psi (\hY(0)),$$
where $f(t) = f(t, X(t), Y(t), Z(t), K(t,\cdot), u(t))$, with $X(t) = X^u(t)$ etc.

By the definition of $H$ we have
\begin{align}\label{eq3.19}
J_1 & = E \left[ \int_0^T \{ H(t) - \hH(t) - \hla(t) (g(t) - \hg(t)) - \hp(t) (b(t) - \hb(t)) \right.\nonumber \\
& \left. \left. - \hq(t) (\sigma(t) - \hsi(t)) - \int_\RB \hr(t,\zeta)(\gamma(t,\zeta) - \hga(t,\zeta)) \nu (d\zeta) \right\} dt \right].
\end{align}
By concavity of $\varphi$, \eqref{eq3.12}, the It\^o formula and \eqref{eq3.18},
\begin{align}\label{eq3.20}
J_2 & \leq E[ \varphi'(\hX(T))(X(T) - \hX(T))] \nonumber \\
& = E [ \hp(T) (X(T) - \hX(T))] - E [ \hla(T) h'(\hX(T))(X(T) - \hX(T))] \nonumber \\
& = \lim_{n \rightarrow \infty} \left( E \left[ \int_0^{\tau_n} \hp(t^-) (dX(t) - d \hX(t)) + \int_0^{\tau_n} (X(t^-) - \hX(t^-)) d \hp (t) \right. \right. \nonumber \\
& + \int_0^{\tau_n} \hq(t) (\sigma(t) - \hsi(t))dt \nonumber \\
& \left. \left. \int_0^{\tau_n} \int_\RB \hr(t,\zeta)(\gamma(t,\zeta) - \hga(t,\zeta)) \nu (d \zeta)\right]\right) - E [ \hla(T) h'(\hX(T))(X(T) - \hX(T))] \nonumber \\
& = E \left[ \int_0^T \hp(t) (b(t) - \hb(t))dt + \int_0^T (X(t)  - \hX(t)) \left( - \frac{\partial \hH}{\partial x}(t)\right)dt \right.\nonumber \\
& \left.\quad + \int_0^T \hq(t)(\sigma(t) - \hsi(t))dt + \int_0^T \int_\RB \hr(t,\zeta) (\gamma(t,\zeta) - \hga(t,\zeta)) \nu(d\zeta) dt \right] \nonumber \\
& \quad - E [ \hla(T) h'(\hX(T))(X(T) - \hX(T))].
\end{align}
By the concavity of $\psi$ and $h$, \eqref{eq3.11} and \eqref{eq3.18},
\begin{align}\label{eq3.21}
J_3 &= \psi(Y(0)) - \psi(\hY(0)) \leq \psi'(\hY(0))(Y(0) - \hY(0)) = \hla(0) (Y(0) - \hY(0)) \nonumber \\
&= \lim_{n \rightarrow \infty}( E[ \hla(\tau_n)(Y(\tau_n) - \hY(\tau_n)) \nonumber \\
& \quad - E \left[ \int_0^{\tau_n} (Y(t^-) - \hY(t^-)) d \hla(t) + \int_0^{\tau_n} \hla(t^-)(dY(t) - d \hY(t)) \nonumber \right.\\
& \quad + \int_0^{\tau_n} \frac{\partial \hH}{\partial z} (t) (Z(t) - \hZ(t))dt \nonumber \\
& \left. \left. + \int_0^{\tau_n} \int_\RB \nabla_k \hH(t,\zeta)(K(t,\zeta) - \hK(t,\zeta)) \nu (d\zeta) dt \right] \right) \nonumber \\
& = E [ \hla(T) (Y(T) - \hY(T))] \nonumber \\
& \quad - E \left[ \int_0^T \frac{\partial \hH}{\partial y}(t) ( Y(t) - \hY(t)) dt + \int_0^T \hla(t)(-g(t) + \hg(t))dt \nonumber \right. \\
& \quad + \int_0^{\tau_n} \frac{\partial \hH}{\partial z}(t) (Z(t) - \hZ(t))dt \nonumber \\
& \quad \left. + \int_0^{\tau_n} \int_\RB \nabla_k \hH(t,\zeta) (K(t,\zeta) - \hK(t,\zeta) \nu(d\zeta) dt \right] \nonumber \\
& \leq E [ \hla(T) h'(\hX(T)) (X(T) - \hX(T))] \nonumber \\
 & \quad - E \left[ \int_0^T \frac{\partial \hH}{\partial y}(t) (Y(t) - \hY(t))dt  \right.
+ \int_0^T \hla(t) (-g(t) + \hg(t))dt \nonumber \\
 & \quad + \int_0^T \frac{\partial \hH}{\partial z}(t) (Z(t) - \hZ(t))dt \nonumber \\
 & \quad \left. + \int_0^T \int_\RB \frac{d \nabla_k \hH}{d \nu}(t,\zeta) (K(t,\zeta) - \hK(t,\zeta)) \nu(d\zeta)dt \right].
\end{align}
Adding \eqref{eq3.19}, \eqref{eq3.20} and \eqref{eq3.21} we get, by \eqref{3.10},
\begin{align}\label{eq3.22}
J(u) & - J(\hu) = J_1 +J_2 + J_3 \nonumber \\
& \leq E \left[ \int_0^T\left\{ H(t) - \hH(t) - \frac{\partial \hH}{\partial x} (X(t) - \hX(t)) \nonumber \right. \right. \\
& \quad - \frac{\partial H}{\partial y}(t) (Y(t) - \hY(t)) - \frac{\partial H}{\partial z}(t) (Z(t) - \hZ(t)) \nonumber \\
& \quad \left.  - \langle \nabla_k \hH(t,\cdot), (K(t,\cdot) - \hK(t,\cdot)\rangle  \biggr{\}} dt \right].
\end{align}
Using that $\hat{\HC}$ is concave, we get by a separating hyperplane argument (see e.g. Rockafellar (1970), Chapt. 5, Sec. 23) that there exists a supergradient $a = (a_0, a_1, a_2, a_3(\cdot)) \in \RB^3 \times \RC$ for $\hat{\HC}(x,y,z,k)$ at $x = \hX(t)$, $y = \hY(t)$, $z = \hZ(t)$ and $k = \hK(t^-,\cdot)$ such that if we define
\begin{align*}
\Phi&(x,y,z,k)  := \hat{\HC}(x,y,z,k) - \hat{\HC}(\hX(t), \hY(t), \hZ(t), \hK(t, \zeta)) \\
 & - [a_0(x-\hX(t)) + a_1(y-\hY(t)) + a_2(z - \hZ(t)) + \int_\RB a_3(\zeta) (k(\zeta) - \hK(t,\zeta))\nu(d\zeta)],
 \end{align*}
 then
 $$\Phi(x,y,z,k) \leq 0 \text{ for all } x,y,z,k.$$
 On the other hand, since
 $$\Phi (\hX(t), \hY(t), \hZ(t), \hK(t,\cdot)) = 0$$
 we get
\begin{align*}
\frac{\partial \hH}{\partial x}(t) & = \frac{\partial \hat{\HC}}{\partial x} (\hX(t), \hY(t), \hZ(t), \hK(t,\cdot)) = a_0 \\
\frac{\partial \hH}{\partial y}(t) & = \frac{\partial \hat{\HC}}{\partial y} (\hX(t), \hY(t), \hZ(t), \hK(t,\cdot)) = a_1 \\
\frac{\partial \hH}{\partial z}(t) & = \frac{\partial \hat{\HC}}{\partial z} (\hX(t), \hY(t), \hZ(t), \hK(t,\cdot)) = a_2 \\
\nabla_k\hH(t,\zeta) & = \nabla_k \hat{\HC} (\hX(t), \hY(t), \hZ(t), \hK(t,\cdot)) = a_3. \\
\end{align*}
If we combine this with \eqref{eq3.22} we obtain 
\begin{align*}
J(u) - J(\hu) & \leq \hat{\HC}(X(t), Y(t), Z(t), K(t,\cdot)) \\
& \quad - \hat{\HC}(\hX(t), \hY(t), \hZ(t), \hK(t,\cdot)) \\
& \quad - \frac{\partial \hat{\HC}}{\partial x} (\hX(t), \hY(t), \hZ(t), \hK(t,\cdot))(X(t) - \hX(t)) \\
& \quad - \frac{\partial \hat{\HC}}{\partial y} (\hX(t), \hY(t), \hZ(t), \hK(t,\cdot))(Y(t) - \hY(t)) \\
& \quad - \frac{\partial \hat{\HC}}{\partial z} (\hX(t), \hY(t), \hZ(t), \hK(t,\cdot))(Z(t) - \hZ(t)) \\
& \quad - \langle \nabla_k \hat{\HC} (\hX(t), \hY(t), \hZ(t), \hK(t,\cdot)),K(t,\cdot) - \hK(t,\cdot)\rangle  \\
& \leq 0, \text{ by concavity of } \HC.
\end{align*}
\fproof

We proceed to a strengthened \emph{necessary} maximum principle. It stronger than the corresponding result in e.g. \O ksendal \& Sulem (2012), because of the weaker growth conditions.

We make the following assumptions:

\begin{itemize}
\item[A1.] For all $t_0 \in [0,T]$ and all bounded $\GC_{t_0}$-measurable random variables $\alpha(\omega)$, the control $\theta(t,\omega) := \chi_{[t_0,T]}(t) \alpha(\omega)$ belongs to $\AC_\GB$.

\item[A2.] For all $u, \beta_0 \in \AC_\GB$ with $\beta_0(t)  \leq K < \infty$ for all $t$, define
 $$\delta(t) := \frac{1}{2 K} dist(u(t), \partial \UB) \wedge 1> 0$$
 and put
 \begin{equation} \label{eq3.24a}
 \beta(t) := \delta(t) \beta_0(t).
 \end{equation}
 Then  the control
$$\tu(t) := u(t) + a \beta(t) \; ; \; t \in [0,T]$$
belongs to $\AC_\GB$ for all $a \in (- 1, 1)$.

\item[A3.] For all $\beta$ as in \eqref{eq3.24a} the derivative processes
\begin{align*}
x(t) & := \frac{d}{da} X^{u+a\beta}(t) \mid_{a=0}, \\
y(t) & := \frac{d}{da} Y^{u+a\beta}(t) \mid_{a=0}, \\
z(t) & := \frac{d}{da} Z^{u+a\beta}(t) \mid_{a=0}, \text{ and } \\
k(t,\zeta) & := \frac{d}{da} K^{u+a\beta}(t,\zeta) \mid_{a=0}
\end{align*}
exists, and belong to $L^2(dm \times dP)$, $L^2(dm \times dP)$, $L^2(dm \times dP)$ and $L^2(dm \times dP \times d\nu)$ respectively, and
\begin{equation}\label{eq3.23}
\begin{cases}
dx(t) =  \displaystyle \left\{ \frac{\partial b}{\partial x}(t) x(t) + \frac{\partial b}{\partial y}(t)y(t) + \frac{\partial b}{\partial z}(t)z(t) + \langle\nabla_k b, k(t,\cdot)\rangle  + \frac{\partial b}{\partial u}(t) \beta(t) \right\} dt \\
 \quad \displaystyle + \left\{ \frac{\partial \sigma}{\partial x}(t) x(t) + \frac{\partial \sigma}{\partial y}(t) y(t) + \frac{\partial \sigma}{\partial z}(t) z(t) + \langle \nabla_k \sigma, k(t,\cdot)\rangle + \frac{\partial \sigma}{\partial u}(t) \beta(t)\right\} dB(t) \\
 \quad \displaystyle + \int_\RB
 \left\{ \frac{\partial \gamma}{\partial x}(t,\zeta) x(t) + \frac{\partial \gamma}{\partial y}(t,\zeta) y(t) + \frac{\partial \gamma}{\partial z}(t,\zeta) z(t) + \langle \nabla_k \gamma(t,\zeta), k(t,\cdot)\rangle \right.\\
 \qquad \left. \displaystyle + \frac{\partial \gamma}{\partial u}(t,\zeta) \beta(t)\right\}  \tN(dt,d\zeta) \; ; \; t \in [0,T] \\
 x(0) = 0
\end{cases}
\end{equation}
\begin{align}\label{eq3.24}
dy(t) & = -\left\{ \frac{\partial g}{\partial x}(t) x(t) + \frac{\partial g}{\partial y}(t)y(t) + \frac{\partial g}{\partial z}(t)z(t) + \langle \nabla_k g(t), k(t,\cdot)\rangle + \frac{\partial g}{\partial u}(t) \beta(t)   \right\} dt \nonumber \\
& + z(t) dB(t) + \int_\RB k(t,\zeta) \tN(dt, d\zeta).
\end{align}
\end{itemize}

\begin{theorem}(Strengthened necessary maximum principle)\label{th3.2}
The following are equivalent:
\begin{myenumerate}
\item $\displaystyle \frac{d}{da} J(u + a \beta) \mid_{a= 0} = 0$ for all bounded $\beta \in \AC_\GB$ of the form \eqref{eq3.24a}.
\item $\displaystyle E \left[ \frac{\partial H}{\partial u}(t) \mid \GC_t \right]= 0$ for all $t \in [0,T]$.
\end{myenumerate}
\end{theorem}

\dproof
Define a sequence of stopping times $\tau_n \; ; \; n = 1, 2, \ldots$ as follows:
\begin{align*}
\tau_n & = \inf \{t > 0 \; ; \; \max \{ |p(t)|, |\nabla \sigma(t)|, \| \nabla \gamma(t,\cdot)\|, |x(t)|, |q(t)|, \\
& \|r(t,\cdot)\|, | \lambda(t)|, |Z(t)|, \|k(t,\cdot)\| \} \geq n \} \wedge T.
\end{align*}
Then it is clear that $\tau_n \rightarrow T$ as $n \rightarrow \infty$
and
\begin{align}\label{eq3.25}
E & \left[ \int_0^{\tau_n} p(t) \sigma(t) dB(t)\right] = E \left[ \int_0^{\tau_n} p(t) \int_\RB \gamma(t,\zeta) \tN (dt,\zeta)  \right] \nonumber \\
& = E\left[ \int_0^{\tau_n} x(t) q(t) dB(t)\right] = 
E \left[ \int_0^{\tau_n}  \int_\RB x(t) r(t,\zeta) \tN (dt,\zeta)  \right] \nonumber \\
&  = E \left[ \int_0^{\tau_n} \lambda(t) Z(t) dB(t) \right] = E \left[ \int_0^{\tau_n} \int_\RB \lambda(t) k(t,\zeta) \tN (dt,d\zeta) \right] = 0 \text{ for all } n.
\end{align}

We can write $\displaystyle \frac{d}{da} J(u + a \beta) \mid_{a=0} = I_1 + I_2 + I_3$, where
\begin{align*}
I_1 & = \frac{d}{da} E \left[ \int_0^T f(t,X^{u+a \beta}(t), Y^{u+a \beta}(t), Z^{u+a \beta}(t), K^{u+a \beta}(t), u(t) + a \beta(t))dt\right]_{a=0} \\
I_2 & = \frac{d}{da} E[ \varphi(X^{u+a \beta}(T))]_{a=0}\\
I_3 & = \frac{d}{da} [ \psi(Y^{u+a \beta}(0))]_{a=0}.
\end{align*}
By our assumptions on $f$, $\varphi$ and $\psi$ we have
\begin{align}\label{eq3.26}
I_1 & = \left[ \int_0^T \left\{ \frac{\partial f}{\partial x}(t) x(t) + \frac{\partial f}{\partial y}(t) y(t) + \frac{\partial f}{\partial z}(t) z(t) + \langle\nabla_k f(t,\cdot), k(t,\cdot)\rangle + \frac{\partial f}{\partial u}(t) \beta(t)\right\} dt\right]\nonumber \\
I_2 & = E[ \varphi' (X(T)) x(T)] = E[p(T) x(T)]  \text{ and } \nonumber \\
I_3 &  = \psi'(Y(0))y(0) = \lambda(0) y(0).
\end{align}
By the It\^o formula and \eqref{eq3.25}
\begin{align}\label{eq3.27}
I_2 & = E[p(T) x(T)] = \lim_{n \rightarrow \infty} E[p(\tau_n) x(\tau_n)] \nonumber \\
& = \lim_{n \rightarrow \infty} E \left[ \int_0^{\tau_n} p(t) dx(t) + \int_0^{\tau_n} x(t) dp(t) + \int_0^{\tau_n} d[p,x](t) \right] \nonumber \\
& = \lim_{n \rightarrow \infty} E \left[ \int_0^{\tau_n} p(t) \left\{ \frac{\partial b}{\partial x}(t) x(t) + \frac{\partial b}{\partial y}(t) y(t) + \frac{\partial b}{\partial z}(t) z(t) + \langle \nabla_k  b(t),k(t,\cdot)\rangle \right.\right.\nonumber \\
& \left. \quad + \frac{\partial b}{\partial u}(t) \beta(t) \right\} dt + \int_0^{\tau_n}x(t) \left( - \frac{\partial H}{\partial x}(t)\right) dt  + \int_0^{\tau_n} q(t) \left\{ \frac{\partial \sigma}{\partial x}(t) x(t)  \right. \nonumber \\
& \quad \left. + \frac{\partial \sigma}{\partial y}(t) y(t) + \frac{\partial \sigma}{\partial z}(t) z(t)  + \langle \nabla_k \sigma(t), k(t,\cdot)\rangle + \frac{\partial \sigma}{\partial u}(t) \beta(t) \right\}dt \nonumber \\
& \quad + \int_0^{\tau_n}\int_\RB r(t,\zeta) \left\{ \frac{\partial \gamma}{\partial x}(t,\zeta)x(t) + \frac{\partial \gamma}{\partial y}(t,\zeta)y(t) + \frac{\partial \gamma}{\partial z}(t,\zeta)z(t) + < \nabla_k \gamma(t,\zeta), k(t,\cdot)> \right.\nonumber \\
& \quad \left. \left. + \frac{\partial \gamma}{\partial u}(t,\zeta)\beta(t)\right\} \nu (d\zeta) dt \right] \nonumber \\
& = \lim_{n \rightarrow \infty} E \left[ \int_0^{\tau_n} x(t) \left\{ \frac{\partial b}{\partial x}(t) p(t) + \frac{\partial \sigma}{\partial x}(t) q(t) + \int_\RB \frac{\partial \gamma}{\partial x}(t,\zeta) r(t,\zeta) \nu(d\zeta) - \frac{\partial H}{\partial x}(t) \right\} dt\right. \nonumber \\
& \quad + \int_0^{\tau_n} y(t) \left\{ \frac{\partial b}{\partial y}(t) p(t) + \frac{\partial \sigma}{\partial y}(t) q(t) + \int_\RB \frac{\partial \gamma}{\partial y}(t,\zeta)  r(t,\zeta) \nu(d\zeta)  \right\} dt\nonumber \\
& \quad + \int_0^{\tau_n} z(t) \left\{ \frac{\partial b}{\partial z}(t) p(t) + \frac{\partial \sigma}{\partial z}(t) q(t) + \int_\RB \frac{\partial \gamma}{\partial z}(t,\zeta)  r(t,\zeta) \nu(d\zeta)  \right\} dt\nonumber \\
& \quad \left. + \int_0^{\tau_n}\int_\RB   \langle k(t,\cdot),  \nabla_kb(t)p(t) + \nabla_k \sigma(t) q(t) + \int_\RB \nabla_k \gamma(t,\zeta) r(t,\zeta) \nu (d \zeta) \rangle \nu(d\zeta)dt\right]\nonumber \\
& = \lim_{n \rightarrow \infty} E \left[ \int_0^{\tau_n} x(t) \left\{ - \frac{\partial f}{\partial x}(t) - \lambda (t) \frac{\partial g}{\partial x}(t)  \right\} dt \right. \nonumber \\
& \quad + \int_0^{\tau_n} y(t) \left\{ \frac{\partial H}{\partial y}(t)- \frac{\partial f}{\partial y}(t) - \lambda(t)  \frac{\partial g}{\partial y}(t)    \right\} dt\nonumber \\
& \quad + \int_0^{\tau_n} z(t) \left\{ \frac{\partial H}{\partial z}(t)- \frac{\partial f}{\partial z}(t) - \lambda(t)  \frac{\partial g}{\partial z}(t)    \right\} dt\nonumber \\
& \quad + \int_0^{\tau_n} \int_\RB k(t,\zeta) \{ \nabla_kH(t) - \nabla_k f(t) - \lambda(t) \nabla_k g(t)\} \nu (d\zeta) dt \nonumber \\
& \quad \left. + \int_0^{\tau_n} \beta(t) \left\{ \frac{\partial H}{\partial u}(t)- \frac{\partial f}{\partial u}(t) - \lambda(t)  \frac{\partial g}{\partial u}(t)    \right\} dt \right]\nonumber \\
& = - I_1 - E \left[ \int_0^T \lambda(t) \left\{ \frac{\partial g}{\partial x}(t) x(t) + \frac{\partial g}{\partial y} (t) y(t) + \frac{\partial g}{\partial z} (t) z(t) \right. \right.\nonumber \\
& \quad \left. + \langle \nabla_kg(t), k(t,\cdot)\rangle + \frac{\partial g}{\partial u}(t) \beta(t) \right\} dt \nonumber \\
& \quad + E \left[ \int_0^T \left\{ \frac{\partial H}{\partial y} (t) y(t) + \frac{\partial H}{\partial z} (t) z(t)  + \langle \nabla_k H(t), k(t,\cdot)\rangle + \frac{\partial H}{\partial u} (t) \beta(t)\right\}dt\right]
\end{align}
Again by the It\^o formula and \eqref{eq3.25},
\begin{align}\label{eq3.28}
I_3 & = \lambda(0) y(0) = \lim_{n \rightarrow \infty} E \left[ \lambda(\tau_n) y(\tau_n)
  -  \left( \int_0^{\tau_n} \lambda(t) dy(t) + \int_0^{\tau_n} y(t) d \lambda(t) + \int_0^{\tau_n} d[\lambda, y](t)\right)\right] \nonumber \\
& = E[ \lambda(T) y(T)] \nonumber \\
&\quad - \lim_{n \rightarrow \infty} \left( E \left[ \int_0^{\tau_n} \lambda(t) \left\{ - \frac{\partial g}{\partial x}(t) x(t) - \frac{\partial g}{\partial y}(t) y(t) - \frac{\partial g}{\partial z}(t) z(t) \right. \right. \right. \nonumber \\
& \left. \quad - \langle \nabla_k g(t), k(t,\cdot)\rangle - \frac{\partial g}{\partial u}(t) \beta(t) \right\}dt \nonumber \\
& \left. \left.\quad + \int_0^{\tau_n} y(t) \frac{\partial H}{\partial y}(t) dt + 
\int_0^{\tau_n} z(t) \frac{\partial H}{\partial z}(t) dt +
\int_0^{\tau_n}   \int_\RB k(t,\zeta) \nabla_k H(t,\zeta) \nu (d \zeta) dt \right] \right).\end{align}
Summing \eqref{eq3.26}, \eqref{eq3.27} and \eqref{eq3.28} we get
$$ \frac{d}{da} J(u + a \beta) \mid_{a=0} = I_1 + I_2 + I_3 = E \left[  \int_0^T \frac{\partial H}{\partial u}(t) \beta(t)dt \right].$$
We conclude that
$$\frac{d}{da} J(u + a \beta) \mid_{a=0}  = 0  $$
if and only if
$$E \left[ \int_0^T \frac{\partial H}{\partial u}(t) \beta(t) dt \right] = 0 \; ; \; \text{ for all bounded } \beta \in \AC_\GB \text { of the form } \eqref{eq3.24a}.$$
In particular, applying this to $\beta(t) = \theta(t)$ as in A1, we get that this is again equivalent to
$$E \left[ \frac{\partial H}{\partial u}(t) \mid \GC_t\right]  = 0 \text{ for all } t \in [0,T].$$
\fproof

\section{Application}\label{sec4}

In this section we give some applications of the theory in Section \ref{sec3}.

\subsection{Utility maximization}\label{sec4.1}

Consider a financial market where the unit price $S_0(t)$ of the risk free asset is
\begin{equation}\label{eq4.a}
S_0(t) = 1 \; ; \; t \in [0,T]
\end{equation}
and the unit price $S_1(t)$ of the risky asset is given by
\begin{equation}\label{eq4.b}
dS_1(t) = S_1(t) [b_0(t) dt + \sigma_0(t) dB(t)] \; ; \; t \in [0,T]
\end{equation}
where $b_0(t), \sigma_0(t)$ are given $\FB$-adapted processes.

Then the wealth process $X_\pi(t)$ associated to a (self-financing) portfolio $\pi(t)$ is given by
\begin{equation}\label{eq4.1}
\begin{cases}
dX_\pi (t) & = \pi(t) X(t^-) [b_0(t) dt + \sigma_0(t) dB(t)] \; ; \; t \geq 0 \\
X_\pi(0) & = x_0 > 0.
\end{cases}
\end{equation}

 Let $U$ be a given utility function. We want to find $\pi^* \in \AC$ such that
\begin{equation}\label{eq4.2}
E[U(X_{\pi^*}(T))] = \sup_{\pi \in \AC} E[U(X_\pi(T))],
\end{equation}
where $\AC$ is the given family of admissible $\FB$-adapted portfolios $\pi$ with values in $\RB$. \\ The Hamiltonian for this problem is
\begin{equation}\label{eq4.3}
H(t,x,y,z,k,\pi, \lambda,p,q,r) = \pi x b_0 p + \pi x \sigma_0(t) q
\end{equation}
and the adjoint equation is
\begin{equation}\label{eq4.4}
\begin{cases}
dp(t) & = - \pi(t) \{b_0(t) p(t) + \sigma_0(t) q(t)\} dt + q(t) dB(t) \; ; \; 0 \leq t \leq T \\
p(T) & = U'(X_\pi(T)).
\end{cases}
\end{equation}
Suppose $\hpi \in \AC$ is an optimal portfolio for the problem \eqref{eq4.2} with corresponding solutions $\hX, \hp, \hq, \hr$ of \eqref{eq4.1} and \eqref{eq4.4}. Then $\displaystyle \frac{\partial \hH}{\partial \pi}(t) = 0$, which gives
\begin{equation}\label{eq4.5} b_0(t) \hp(t) + \sigma_0(t) \hq(t) = 0.
\end{equation}
Hence
$$\hq(t) = - \frac{b_0(t)}{\sigma_0(t)} \hp(t)$$
and \eqref{eq4.4} becomes
\begin{equation}\label{eq4.6}
\begin{cases}
d\hp(t) & \displaystyle = \hq(t) dB(t) = - \frac{b_0(t)}{\sigma_0(t)} \hp(t) dB(t) \; ; \; 0 \leq t \leq T \\
\hp(T) & = U'(\hX(T)).
\end{cases}
\end{equation}
Define
\begin{equation}\label{eq4.7}
\theta_0(t) = - \frac{b_0(t)}{\sigma_0(t)}.
\end{equation}
Then by \eqref{eq4.5}
$$b_0(t) + \sigma_0(t) \theta_0(t) = 0,$$
and the Girsanov theorem gives that the measure $Q$ defined by
\begin{equation}\label{eq4.8}
dQ = \Gamma(T) dP \text{ on } \FC_T
\end{equation}
is an equivalent local martingale measure, where $\Gamma(t) = \Gamma_{\theta_0}(t)$ is given by
\begin{equation}\label{eq4.9}
d \Gamma(t)  = \Gamma(t) \theta_0(t) dB(t) \; ; \; 
\Gamma(0)  = 1.
\end{equation}
Note that
$$\frac{d \hp(t)}{\hp(t)} = \frac{d \Gamma(t)}{\Gamma(t)}$$
so
$$\hp(t) = E[U'(\hX(T))] \Gamma(t).$$
By \eqref{eq4.7} and \eqref{eq4.9} we have
\begin{equation}\label{eq4.10}
\Gamma(t) = \exp \left( - \int_0^t \frac{b_0(s)}{\sigma_0(s)} dB(s) - \frac{1}{2} \int_0^t \frac{b_0^2(s)}{\sigma_0^2(s)} ds \right) \; ; \; 0 \leq t \leq T.
\end{equation}
Hence
$$\hp(T) = U'(\hX(T)) = E[U'(\hX(T))] \Gamma(T),$$
i.e.
\begin{equation}\label{eq4.10a}\hX(T) = I(c \;  \Gamma(T))\end{equation}
where
$$I = (U')^{-1} \text{ and } c = E[U'(\hX(T))].$$
It remains to find $c$. We can write \eqref{eq4.1} as
$$
\begin{cases}
d \hX(t) & = \hpi(t) \hX(t) [b_0(t) dt + \sigma_0(t) dB(t)] \; ; \; 0 \leq t \leq T \\
\hX(T) & = I(c \Gamma(T)).
\end{cases}
$$
If we define
\begin{equation}\label{eq4.11}
\hZ(t) = \hpi(t) \hX(t) \sigma_0(t)
\end{equation}
this becomes a BSDE
\begin{equation}\label{eq4.12}
\begin{cases}
d \hX(t) & \displaystyle = \frac{\hZ(t) b_0(t)}{\sigma_0(t)} dt + \hZ(t) db(t) \; ; \; 0 \leq t \leq T \\
\hX(T) & = I(c \Gamma(T)).
\end{cases}
\end{equation}
The solution of this BSDE is
\begin{equation}\label{eq4.13}
\hX(t) = \frac{1}{\Gamma(t)} E[I(c \Gamma(T)) \Gamma(T) \mid \FC_t].
\end{equation}
In particular,
\begin{equation}\label{eq4.14}
x = \hX(0) = E[I(c \Gamma(T)) \Gamma(T)].
\end{equation}
This is an equation which (implicitly) determines the value of $c$.
When $c$ is found, we have the optimal terminal wealth $\hX(T)$ given by \eqref{eq4.10a}. Solving the resulting BSDE for $\hZ(t)$, we get the corresponding optimal portfolio $\hpi(t)$ by \eqref{eq4.11}.
We have proved:
\begin{theorem}\label{th4.1}
The optimal terminal wealth $\hX(T) = X_{\pi^*}(T)$ for the portfolio optimization problem \eqref{eq4.2} is given by \eqref{eq4.10a}, where the constant $c > 0$ is found implicitly by equation \eqref{eq4.14}.
\end{theorem}

\subsection{Risk minimization}\label{sec4.2}
The starting point is the same as in Section \ref{sec4.1}, with a wealth equation given by \eqref{eq4.1}.
This time we want to \emph{minimize the risk} $\rho(X_\pi(T))$ of the terminal value $X_\pi(T)$, defined by
\begin{equation}\label{eq4.15}
\rho(X_\pi(T)) = - Y_\pi(0)
\end{equation}
where
\begin{equation}\label{eq4.16}
\begin{cases}
dY_\pi(t) & = - g(Z(t))dt + Z(t) dB(t) \; ; \; t \in [0,T] \\
Y_\pi(T) & =  X_\pi(T),
\end{cases}
\end{equation}
for some given concave function $g$.
Thus we want to find $\hpi \in \AC$ and $\rho(X_{\hpi}(T)) := - Y_{\hpi}(0)$ such that
\begin{equation}\label{eq4.16a}
\inf_{\pi \in \AC} (- Y_\pi(0)) = - Y_{\hpi}(0).
\end{equation}
In this case the Hamiltonian becomes
\begin{equation}\label{eq4.17}
H(t,x,y,z,k,\pi,\lambda,p,q,r) = \pi x b_0(t) p + \pi x \sigma_0(t) q + \lambda g(z).
\end{equation}
The adjoint equations are (see \eqref{eq3.11} - \eqref{eq3.12})
\begin{equation}\label{eq4.18}
\begin{cases}
dp(t) & = - \{ \pi(t) b_0(t) p(t) + \pi(t) \sigma_0(t) q(t)\} dt + q(t)  dB(t) \\
p(T) & =  \lambda (T)
\end{cases}
\end{equation}
and
\begin{equation}\label{eq4.19}
\begin{cases}
d \lambda (t) & = \lambda(t) g' (Z(t)) dB(t) \\
\lambda (0) & = 1
\end{cases}
\end{equation}
i.e.
\begin{equation}\label{eq4.20}
\lambda(t) = \exp \left( \int_0^t g'(Z(s)) dB(s) - \frac{1}{2} g' (Z(s))^2ds \right).
\end{equation}
If $\hpi$ is optimal, then
\begin{equation}\label{eq4.21}
b_0(t) \hp(t) + \sigma_0(t) \hq(t) = 0.
\end{equation}
This gives
\begin{equation}\label{eq4.22}
\begin{cases}
d\hp(t) & \displaystyle = \hq(t) dB(t) = - \frac{b_0(t)}{\sigma_0(t)} \hp(t) dB(t) \; ; \; 0 \leq t \leq T \\
\hp(T) & =  \hla(T).
\end{cases}
\end{equation}
Comparing with \eqref{eq4.19} we see that the solution $(\hp, \hq)$ of the BSDE \eqref{eq4.22} is
\begin{equation}\label{eq4.23}
\hp(t) =  \hla(t), \hq(t) =  \hla(t) g' (\hZ(t)).
\end{equation}
Hence by \eqref{eq4.21}
\begin{equation}\label{eq4.24}
g'(\hZ(t)) = - \frac{b_0(t)}{\sigma_0(t)}.
\end{equation}
If, for example, 
\begin{equation}\label{eq4.25}
g(z) = - \frac{1}{2} z^2
\end{equation}
then \eqref{eq4.24} gives
$$\hZ(t) =  \frac{b_0(t)}{\sigma_0(t)}.$$
Substituted into \eqref{eq4.1} this gives, using \eqref{eq4.10} (with $\Gamma(t)$ as in \eqref{eq4.10}),
\begin{align}\label{eq4.26}
 \hX(T) = \hY(T) & = \hY(0) + \int_0^T  \frac{1}{2} \left( \frac{b_0(s)}{\sigma_0(s)}\right)^2 ds + \int_0^T \frac{b_0(s)}{\sigma_0(s)} dB(s) \nonumber \\
 &= \hY(0) - \ln \Gamma(T).
 \end{align}
 We take expectation w.r.t. the martingale measure
 \begin{equation}\label{eq4.27}
 dQ(\omega) = \Gamma(T) dP
 \end{equation}
 as in \eqref{eq4.8} and get
 \begin{equation}\label{eq4.28}
 -\hY(0) = - x - E_Q [ \ln \Gamma(T)] = - x - E \left[ \frac{dQ}{dP} \ln \frac{dQ}{dP}\right].
 \end{equation}
 Note that $\displaystyle H(Q \mid P) := E \left[ \frac{dQ}{dP} \ln \frac{dQ}{dP}\right]$ is the {\it entropy} of $Q$ with respect to $P$.
 
 Now that the optimal valuee $\hY(0)$ is found, we can use \eqref{eq4.26} to find the corresponding optimal terminal wealth $\hX(T)$, and from there the optimal portfolio as we did in Example \ref{sec4.1}. We have proved:

 \begin{theorem}\label{th4.2}
 
 Suppose \eqref{eq4.25} holds. Then the minimal risk $-Y_{\hpi}(0) = - \hY(0)$ of problem \eqref{eq4.16a} is given by \eqref{eq4.28}, where $dQ = \Gamma(T) dP$ is the unique equivalent martingale measure for the market \eqref{eq4.a}-\eqref{eq4.b}.
 \end{theorem}
 
\subsection{The newsvendor problem}\label{sec4.3}
 
 Consider the following Stackelberg game (see  [\O ksendal, Sandal \& Ubøe (2013)], in which the two players, the {\it leader} and the {\it follower}, are the following:
 
 \begin{myenumerate}
 \item  The manufacturer (leader), who decides the wholesale price $w_t$ per unit.
 \item The retailer (follower), who decides the quantity to order, $Q_t$, and the retail price $R_t$, for delivery $\delta > 0$ (fixed) units of time later.
 \end{myenumerate}
 
 The demand process $X_t$ is assumed to satisfy
 \begin{equation}\label{eq4.29}
 dX_t = (K - R_t) dt + \sigma dB_t.
 \end{equation}
 We introduce the following quantities:
 \noindent $M$ is the  production cost per unit (fixed), and 
 \noindent $S$ is the  salvage price per unit (fixed).
 
 For given wholesale (leader L) price $w_t$ the retailer (follower F) tries to maximize
 \begin{equation}\label{eq4.30}
 J_F(w,Q,R) := E \left[ \int_\delta^T \{(R_t-S) \min [X_t, Q_t] - (w_t - S) Q_t \} dt \right]
 \end{equation}
 over all $\GC_t := \FC_{t- \delta}$ -adapted processes $Q_t, R_t$.
 Let $\hQ_t = \hQ(w)_t$, $\hR_t = \hR(w)_t$ be the optimizers of $J_F(w,Q,R)$. The leader knows that this will be the response from the (rational) follower to her choice of $w$, so she tries to find $\hw_t$ which maximizes
 \begin{equation}\label{eq4.31}
 w \rightarrow J_L(w, \hQ(w), \hR(w)) := E \left[ \int_\delta^T (w_t - M) \hQ(w)_t dt \right].
 \end{equation}
 How do we find $\hQ, \hR$ and then $\hw$ ?
 
 We first consider the follower problem, which is for given $w_t$ to find processes $\hQ, \hR$ such that
 \begin{equation}\label{eq4.32}
 \sup_{Q_t, R_t \in \AC_\GB} J_F (w,Q,R) = J_F (w, \hQ, \hR).
 \end{equation}
 This is an SDE stochastic control problem  of classic type.
 The Hamiltonian is
 $$H_F(t,x,Q,R,p,q) = (R-S) \min [x,Q] - (w_t - S) Q + (K-R)p + \sigma q$$
 and the adjoint equation is
 $$\begin{cases}
 dp_F(t) & = - (R_t - S) \chi_{[0, Q_t]}(X_t) dt + q_F(t) dB_t \\
 p_F(T) & = 0,
 \end{cases}
 $$
 which has the solution
 \begin{equation}\label{eq4.33}
 p_F(t) = E \left[ \int_t^T (R_s - S) \chi_{[0,Q_S]}(X_s) ds \mid \FC_t \right].
 \end{equation}
 The first order conditions for a maximum $Q_t, R_t$ of $H$ are
 \begin{equation}\label{eq4.34}
 E \left[ \chi_{[0,X_t^+]} (Q_t) \mid \GC_t\right] = \frac{R_t - w_t}{R_t-S} = \frac{R_t-w_t}{(R_t-w_t) + w_t-S}
 \end{equation}
 and
 \begin{equation}\label{eq4.35}
 E\left[ \min [X_t, Q_t] - \int_t^T (R_s - X_s) \chi_{[0,Q_s]} (X_s)ds \mid \GC_t \right] = 0.
 \end{equation}
 These equations can be solved as follows : Define
 \begin{equation}\label{eq4.36}
 h_t(x) = E\left[\chi_{[0,x]}(X_t) \mid \GC_t\right]
 \end{equation}
 and
 \begin{equation}\label{eq4.37}
 f_t(x) = E\left[X_t\chi_{[0,x]}(X_t) \mid \GC_t\right]
 \end{equation}
Then
\begin{equation}\label{eq4.38}
\hQ(t) = \hQ_t(w_t) = h^{-1}_t \left( \frac{y}{y+ w_t-S}\right)_{y = \hR_t - w_t}
 \end{equation}
and
\begin{equation}\label{eq4.39}
\hR_t = w_t + F^{-1}_t (w_t, Y_t), \text{ where } Y_t = E [ \min (X_t, Q_t) \mid \GC_t]
 \end{equation}
and
\begin{equation}\label{eq4.40}
F_t = F_t^{(w)}(y) = h^{-1}_t \left( \frac{y}{y+ w_t-S}\right) \frac{w-S}{y+w-S} + f_t \left( h^{-1}_t\left( \frac{y}{y+ w_t-S}\right)\right)
 \end{equation}
 and by \eqref{eq4.35} $Y_t$ is the solution of the BSDE
 \begin{equation}\label{eq4.41}
 Y_t = E \left[ \int_t^T F^{-1}_s (w_s, Y_s) ds \mid \GC_t\right]
 \end{equation}
 i.e.
 \begin{equation}\label{eq4.42}
 dY_t  = - F^{-1}_s (w_t, Y_t) dt + d M_t \; ; \; 
 Y_T  = 0
 \end{equation}
where $M_t$ is an $\GB$-martingale.
Thus
\begin{equation}\label{eq4.42a}
\hQ_t = \hQ_t(\hw_t) = h^{-1}_t \left( \frac{F^{-1}_t (w_t,y)}{F^{-1}_t (w_{t},y) + w_t - S} \right)_{y = Y_t}.
\end{equation}
Therefore, if the BSDE \eqref{eq4.34} is solved for $Y$, we have $\hR_t$ by \eqref{eq4.39} and $\hQ_t$ by \eqref{eq4.38}.

This gives the optimal response $\hQ_t = \hQ(w)_t$ and $\hR_t = \hR(w)_t$ of the follower to the manufacturer's price process $w(\cdot)$.

It remains to find the optimal $\hw$ of the manufacturer. Thus we want to find $\hw \in \AC_L(\EC)$ such that
\begin{equation}\label{eq4.43}
\sup_{w \in \AC_L(\EC)} J(w) = J(\hw)
\end{equation}
where
\begin{equation}\label{eq4.44}
J(w) = E \left[ \int_\delta^T (w_t - M) \hQ(w)_t dt \right].
\end{equation}
Substituting for $\hQ(w)$ we see that this can be written
\begin{equation}\label{eq4.45}
J(w) =  E \left[ \int_\delta^T (w_t - M) h^{-1}_t \left( \frac{F^{-1}_t (w_t, Y_t)}{F^{-1}_t (w_t, Y_t) + w_t-S}\right)dt \right].
\end{equation}
Accordingly, the problem to maximize $J(w)$ is the problem of optimal control of the following {\it fully coupled} system of FBSDEs:
$$\text{(Forward )} \begin{cases}
dX_t & = (K - w_t - F^{-1}_t (w_t, Y_t))dt + \sigma dB_t \; ; \; t \geq 0 \\
X_0 & = x > 0
\end{cases}
$$
$$\text{(Backward )} \begin{cases}
dY_t & =  - F^{-1}_t (w_t, Y_t)dt + dM_t \; ; \; 0 \leq t \leq T \\
Y_T & =  0.
\end{cases}
$$
To handle this we use the method from Section \ref{sec3}:
Define the Hamiltonian by
\begin{align}\label{eq4.46}
H(t,x,y,z,w,\lambda,p,q)& = (w - M) h^{-1}_t \left( \frac{F^{-1}_t (w_t, y)}{F^{-1}_t (w_t, y) + w_t-S}\right) + \lambda F^{-1}_t (w_t,y) \nonumber \\
 &\quad + (K - w - F^{-1}_t (w,y)) p + \sigma q.
\end{align}
The adjoint process $\lambda(t)$ is given by
\begin{equation}\label{4.47}
\begin{cases}
d \lambda(t) & \displaystyle = \frac{\partial H}{\partial y}(t) dt + \frac{\partial H}{\partial z}(t) dB_t = \frac{\partial H}{\partial y}(t) dt \\
\lambda(0) & = 0
\end{cases}
\end{equation}
and the adjoint processes $p(t), q(t)$ are given by
$$\begin{cases}
dp(t) &\displaystyle  = - \frac{\partial H}{\partial x}(t) dt + q(t) dB(t) = q(t) dB(t) \; ; \; 0 \leq t \leq T \\
p(T) & = 0
\end{cases}
$$
which gives $p=q=0$.

Therefore the first order condition for a maximum $\hw$ of the Hamiltonian \eqref{eq4.46} is
\begin{equation}\label{eq4.48}
(\hw - M) \frac{d}{dw}\hQ_t(w)_{w = \hw_t} + \hQ_t (\hw_t ) + \hla(t) \frac{d}{dw} (F^{-1}_t(w, Y_t))_{w=\hw_t} = 0.
\end{equation}
The solution of this equation gives us $\hw_t$ as a function of $Y_t$.
Summarizing, we have proved the following (see \cite{OSU} for details):
\begin{theorem}\label{th4.3}
Suppose the Stackelberg/newsvendor problem \eqref{eq4.30}-\eqref{eq4.31} has a solution. Then the optimal manufacturer's price $\hw_t$ is a solution of equation \eqref{eq4.48}, and the corresponding optimal responses $\hQ, \hR_t$ of the retailer are given vy \eqref{eq4.42a} and \eqref{eq4.39}, respectively, where $Y$ is the solution of the BSDE \eqref{eq4.42} with $w_t = \hw_t$.
\end{theorem}

\subsection{Maximizing the recursive utility}\label{sec4.4}

This example is taken from a paper by K. Aase (2012).
In the following $f : \RB^4  \rightarrow \RB$ is a given function and $\pi_t, c_t$ are given adapted processes.
The representative agent's problem is to solve
\begin{equation}\label{eq4.49}
\begin{cases}
\sup_{\tc} Y^{\tc}(0), \text{ subject to } \\
\displaystyle E \left[\int_0^T \pi_t \tc_t dt \right] \leq E \left[ \int_0^T \pi_t c_tdt \right] \; \text{ (the budget constraint)}.
\end{cases}
\end{equation}
Hence $Y(t) = Y^{\tc}(t)$ is the solution of the BSDE
\begin{equation}\label{4.50}
\begin{cases}
dY(t) & = - f(t, \tc_t, Y(t), Z(t), K(t, \cdot))dt \\
 & \displaystyle + Z(t) dB(t) + \int_{\RB_0} K(t,\zeta) \tN (dt,d\zeta) \; ; \; 0 \leq t \leq T \\
 Y(T) & = 0.
\end{cases}
\end{equation}
This is a constrained recursive utility optimization problem.
For $\ell > 0$ define the Lagrangian
\begin{equation}\label{eq4.51}
\LC_\ell(\tc) = Y^{\tc}(0) - \ell E \left[\int_0^T \pi_t (\tc_t - c_t) dt \right].
\end{equation}
Suppose we for each $\ell > 0$ can find an optimal $\tc_t(\ell)$ such that
\begin{equation}\label{eq4.52}
\sup_{\tc} \LC_\ell(\tc) = \LC_\ell (\tc(\ell))
\end{equation}
without constraints.
Next, suppose we can find $\hel$ such that
\begin{equation}\label{eq4.53}
E \left[\int_0^T \pi_t (\tc_t(\hel) - c_t) dt \right] = 0.
\end{equation}
Then
$$\tc^* := \tc(\hel)$$
is optimal for the original constrained problem. To see this, note that for all $\tc$ we have
\begin{align*}
Y^{\tc(\hel)}(0) & = Y^{\tc(\hel)}(0) - \hel E \left[\int_0^T \pi_t (\tc_t(\hel) - c_t) dt \right] \\
 & = \LC_{\hel}(\tc(\hel)) \geq \LC_{\hel}(\tc)\\
 & = Y^{\tc}(0) - \hel E \left[\int_0^T \pi_t (\tc_t - c_t) dt \right] \\
  & \geq Y^{\tc}(0).
\end{align*}
In view of this, we can solve the original constrained problem \eqref{eq4.49} in two steps:
\begin{description}
\item{\it Step 1.} Maximize $\LC_\ell(\tc)$ over all $\tc$ (without constraints), for each given $\ell > 0$. Call the maximum $\tc(\ell)$.
\item{\it Step 2.} Find $\hel$ such that $\displaystyle E \left[\int_0^T \pi_t (\tc_t(\hel) - c_t) dt \right] = 0$. Then $c^* := \tc(\hel)$ solves the original constrained problem.
\end{description}

We now apply this to the problem \eqref{eq4.49}. Thus we fix $\ell$ and proceed to study the unconstrained problem \eqref{eq4.52} in the context of Section \ref{sec3}:
Consider the controlled FBSDE system consisting of
\begin{equation}\label{eq4.54}
dX(t) = 0 \; ; \; X(0) = 0
\end{equation}
\begin{equation}\label{eq4.55}
\begin{cases}
dY(t) &= - f(t, \tc_t, Y(t), Z(t), K(t,\cdot))dt \\
& \displaystyle Z(t) dB(t) + \int_\RB K(t,\zeta) \tN (dt, d\zeta) \; ; \; 0 \leq t \leq T \\
Y(T) & = 0.
\end{cases}
\end{equation}
The performance functional is
\begin{equation}\label{eq4.56}
J(\tc) = - \ell E \left[ \int_0^T  \pi_t (\tc(t) - c(t))dt \right] + Y^{\tc}(0),
\end{equation}
where $\ell > 0$ is the Lagrange multiplier.
The Hamiltonian for this problem is
\begin{equation}\label{eq4.57}
H(t,y,z,k,\lambda,p,p,r) = - \ell \pi(t) (\tc - c(t)) + \lambda f(t, \tc,y,z,k)
\end{equation}
and the adjoint equation is
\begin{equation}\label{eq4.58}
\begin{cases}
d\lambda(t) &  = \displaystyle\lambda(t) \left( \frac{\partial f}{\partial y} (t, \tc(t), Y(t), Z(t),K(t,\cdot)) dt + \frac{\partial f}{\partial z} (t, \tc(t), Y(t), Z(t), K(t,\cdot))dB(t)\right. \\
& \qquad \displaystyle \left. + \int_\RB \frac{d \nabla_k f}{d \nu}(t, \tc(t), Y(t), Z(t), K(t,\cdot)) (\zeta) \tN(dt, d\zeta)\right) \; ; \; 0 \leq t \leq T \\
\lambda(0) & = 1.
\end{cases}
\end{equation}

As before $\nabla_kf$ denotes the Fr\'echet derivative of $f$ with respect to $k : \RB_0 \rightarrow \RB$ and $\displaystyle \frac{d \nabla_k f}{d \nu}(\zeta)$ denotes its Radon-Nikodym derivative with respect to $\nu$.

If $\hc = \hc(\ell)$ is optimal for a given Lagrange multiplier $\ell$, with corresponding values $\hla(t), \hY(t), \hZ(t), \hK(t,\cdot)$ we get by the It\^o formula that the solution of \eqref{eq4.58} is
\begin{align}\label{eq4.59}
\hla(t) & = \exp \left( \int_0^t \left\{ \frac{\partial \hf}{\partial y}(s) - \frac{1}{2} \left( \frac{\partial \hf}{\partial z}\right)^2 + \int_{\RB_0}\left( ln \left( 1 + \frac{d \nabla_k \hf}{d \nu}(s,\zeta)\right) - \frac{d \nabla_k \hf}{d \nu}(s,\zeta)\right) \nu (d \zeta) \right\} ds \right. \nonumber \\
& +\left. \int_0^t \frac{\partial \hf}{\partial z}(s) dB(s) + \int_0^t \int_\RB ln\left(1 + \frac{d \nabla_k \hf}{d \nu} (s,\zeta)\right) \tN (ds, d\zeta)\right) \; ; \; t \geq 0,
\end{align}
where we have used the simplified notation
$$\frac{\partial \hf}{\partial y}(s) = \frac{\partial f}{\partial y}(s, \hc(s), \hY(s), \hZ(s), \hK(s, \cdot)) \text{ etc.}$$
and we assume that
\begin{equation}\label{eq4.60}
\frac{d \nabla_k \hf}{d \nu} (s,\zeta) >  - 1 \text{ for all } s,\zeta \text{ a;s.}.
\end{equation}
Maximizing the Hamiltonian with respect to $\tc$ gives the first order equation
$$ - \ell \pi_t + \hla(t) \frac{\partial f}{\partial \tc}(t, \tc(t), \hY(t), \hZ(t), \hK(t,\cdot)) = 0$$
or
\begin{equation}\label{eq4.61}
\ell \pi_t = \hla(t) \frac{\partial f}{\partial \tc}(t, \tc(t), \hY(t), \hZ(t), \hK(t, \cdot)).
\end{equation}
Assume that for each $t,y$ and $z$ the function
$$\tc \rightarrow \frac{\partial f}{\partial c} (t, \tc,y,z,k)$$
has an inverse, denoted by $\displaystyle \left( \frac{\partial f}{\partial \tc}\right)^{-1}(t, \cdot, y,z)$. Then the solution $\tc(t)$ of the first order condition \eqref{eq4.61} can be written
\begin{align}\label{eq4.62}
\tc(t) & = \tc(t, \hla, \hY(t), \hZ(t), \hK(t,\cdot)) \nonumber \\
 & = \left(\frac{\partial f}{\partial \tc}\right)^{-1} \left(t, \frac{\ell \pi_t}{\hla(t)}, \hY(t), \hZ(t), \hK(t,\cdot)\right) =: \tc(t).
 \end{align}
 Substituting this into the forward equation \eqref{eq4.58} for $\hla(t)$ we get
\begin{equation}\label{eq4.63}
\begin{cases}
d\hla(t) &  = \displaystyle\hla(t) \left( \frac{\partial \hf}{\partial y} (t, \tc(t)) dt + \frac{\partial \hf}{\partial z} (t, \tc(t))dB(t)\right. \\
& \qquad \displaystyle \left. + \int_\RB \frac{d \nabla_k \hf}{d \nu}(t, \tc(t)) (\zeta) \tN(dt, d\zeta)\right) \; ; \; 0 \leq t \leq T \\
\lambda(0) & = 1
\end{cases}
\end{equation}
where $\displaystyle  \frac{\partial \hf}{\partial y}(t, \hc(t)) =  \frac{\partial f}{\partial y}(t, \tc(t), \hla(t), \hY(t), \hZ(t), \hK(t,\cdot)), \hY(t), \hZ(t), \hK(t,\cdot)$) etc. This is coupled to the following BSDE for $(\hY, \hZ, \hK)$:
\begin{equation}\label{eq4.64}
\begin{cases}
d\hY(t) & = - f(t, \tc(t, \hla(t), \hY(t), \hZ(t)), \hK(t,\cdot), \hY(t), \hZ(t), \hK(t,\cdot)) dt \\
& \displaystyle + \int_\RB \hZ(t) dB(t) + \int_\RB \hK(t,\zeta) \tN (dt, d\zeta) \; ; \; 0 \leq t \leq T \\
\hY(T) = & 0.
\end{cases}
\end{equation}
We see that $(\hla(t), \hY(t), \hZ(t))$ can be found as the solution of the coupled system of the FBSDE \eqref{eq4.63} and \eqref{eq4.64}. There is a general theory for the solution of coupled systems of FBSDEs. See e.g. Hu \& Peng (1995).

As an important special case, assume that $f$ has the form
\begin{equation}\label{eq4.65}
f(t, \tc,y,z,k) = f_0(t, \tc,y) - \frac{1}{2} A(y) Z^2 - \frac{1}{2} \int_{\RB_0} A_1(y,\zeta) k^2(\zeta) \nu (d\zeta)
\end{equation}
where $f_0$ does not depend on $z$ and $k$.
In this case we see that $\nabla_kf$ is the linear operator
$$h \rightarrow (\nabla_k f) (h) = - \int_\RB A_1(y,\zeta) k(\zeta) h(\zeta) \nu (d\zeta) \; ; \; h \in L^2(\nu).$$
Therefore, as a random measure we have
$$\nabla_k f << \nu,$$
with Radon-Nikodym derivative
\begin{equation}\label{eq4.66}
\frac{d \nabla_k f}{d \nu} (\zeta) = - A_1(y, \zeta) k(\zeta).
\end{equation}

Hence \eqref{eq4.63} gets the form
\begin{equation}\label{eq4.67}
\begin{cases}
d \hla(t) & = \hla(t)\left( \left\{ \frac{\partial f_0}{\partial y}(t, \tc, \hY(t))- \frac{1}{2} A'(\hY(t)) \hZ^2(t) - \frac{1}{2} \int_\RB A'_1 (\hY(t,\zeta) \hK^2(t,\zeta) \nu (d \zeta)\right\} dt \right. \\
&\displaystyle  \left. - A (\hY(t)) \hZ(t) dB(t) - \int_\RB A_1(\hY(t),\zeta) \hK(t,\zeta) \tN (dt, d \zeta)\right) \; ; \; t \geq 0 \\
\hla(0) & = 1
\end{cases}
\end{equation}
which has the solution
\begin{align}\label{eq4.68}
\hla(t) & = \exp \left( \int_0^t \left\{ - \frac{1}{2} (A'(\hY(t)) + A^2(\hY(t))) \hZ^2(t) \right.\right.\nonumber \\
& - \frac{1}{2} \int_\RB A'_1 (\hY(t), \zeta) \hK^2 (t,\zeta) \nu (d \zeta) \nonumber \\
& \left.+ \int_\RB \{ ln (1 - A_1 (\hY(t), \zeta) \hK(t,\zeta)) + A_1(\hY(t),\zeta) \hK(t,\zeta)\} \nu d\zeta) \right\} dt \nonumber \\
& \left. - \int_0^t A(\hY(t)) \hZ(t) dB(t) + \int_0^t \int_\RB ln (1 - A_1 (\hY(t),\zeta) \hK(t,\zeta)) \tN (dt,d\zeta)\right) \; ; \; t \geq 0.
\end{align}
The corresponding BSDE \eqref{eq4.64} for $(\hY, \hZ, \hK)$ gets the form
\begin{equation}\label{eq4.69}
\begin{cases}
d\hY(t) & \displaystyle = - \left\{ f_0(t, \tc(t), \hY(t)) - \frac{1}{2} A (\hY(t)) \hZ^2(t) \right. \\
& \left. \displaystyle - \int_\RB \frac{1}{2} A_1 (\hY(t), \zeta) \hK^2(t,\zeta) \nu (d\zeta) \right\} dt \\
& \displaystyle + \hZ(t) dB(t) + \int_\RB \hK(t,\zeta) \tN (dt, d\zeta) \; ; \; 0 \leq t \leq T \\
\hY(T) & = 0,
\end{cases}
\end{equation}
where
$$\tc(t) = \tc(t, \hla(t), \hY(t), \hZ(t), \hK(t,\cdot)).$$
\section{Risk minimization and stochastic differential games}\label{sec6}

\subsection{A dual representation of convex risk measures}\label{sec6.1}

In Section \ref{sec2} we saw that BSDEs can be used to define convex risk measures. With this representation the problem of risk minimization becomes a problem of stochastic control of FBSDE's, as we saw in Sections \ref{sec3} and \ref{sec4}. 
There is another representation of convex risk measures based on convex duality methods. It goes as follows:
\begin{theorem}[Convex duality representation of convex risk measures]\label{th6.1}
Every convex risk measure $\rho : L^\infty(\FC_T, P) \rightarrow \RB$ is of the following form:
\begin{equation}\label{eq6.1}
\rho(X) = \sup_{Q \in \PC} \{E_Q[-X] - \alpha(Q)\} \; ; \; X \in L^\infty(\FC_T,P)
\end{equation}
where $\PC$ is a family of probability measures $Q \ll P$ and $\alpha : \PC \rightarrow \RB$ is a convex function, usually called a {\it penalty function}. See F\"ollmer \& Schied (2002), Frittelli \& Rosazza-Gianin (2002).
\end{theorem}

For example, $\PC$ could be the set $\PC_\Theta$ of probability measures $Q_\theta$ defined by
\begin{equation}\label{eq6.2}
dQ_\theta(\omega) = M_\theta(T) dP(\omega) \text{ on } \FC_T,
\end{equation}
where $\theta$ denotes the process $(\theta_0(t), \theta_1(t,\zeta))$ and $M_\theta(t)$ is the martingale defined by
\begin{equation}\label{eq6.3}
\begin{cases}
dM_\theta(t) & \displaystyle = M_\theta(t^-) \left[ \theta_0(t) dB(t) + \int_\RB \theta_1(t,\zeta) \tN (dt,d\zeta)\right] \; ; \; t \geq 0 \\
M_\theta(0) & = 1.
\end{cases}
\end{equation}
Here $\Theta$ is the set of admissible scenarios, defined by 
\begin{align}\label{eq6.4}
\Theta = \{ \theta(t) = (\theta_0, \theta_1(t,\zeta)); \; \theta \text{ is $\FC_t$-predictable, $\theta_1(t,\zeta) > -1$ and} \nonumber \\
E \left. \left[ \int_0^T \{ \theta^2_0(t) + \int_\RB \theta^2_1(t,\zeta) \nu(d\zeta)\} dt \right] < \infty \right \}.
\end{align}
A popular choice of the penalty function $\alpha$ is the {\it entropic penalty} $\alpha_e$, defined by
\begin{equation}\label{eq6.5}
\alpha_e(Q) = E \left[ \frac{dQ}{dP} \log \frac{dQ}{dP} \right], \quad Q \in \PC.
\end{equation}
In information theory the quantity $\displaystyle E \left[ \frac{dQ}{dP} \log \frac{dQ}{dP} \right]$ is known as the {\it entropy} of $Q$ with respect to $P$ and it is denoted by $H(Q\mid P)$. The corresponding convex risk measure in \eqref{eq6.1}, with $\PC$ and $\alpha_e$ as above, is called the {\it entropic risk measure} and denoted by $\rho_e$.

We will not give the proof  of the representation \eqref{eq6.1} here, but we refer to F\"ollmer \& Schied (2002), Frittelli \& Rosazza-Gianin (2002) for details. See also F\"ollmer \& Schied (2010).
Note that if $\alpha = 0$ then the corresponding risk measure will be {\it sub-additive}, i.e.
\begin{equation}\label{eq6.6}
\rho (X + Y) \leq \rho (X) + \rho(Y)
\end{equation}
and {\it positive homogeneous}, i.e.
\begin{equation}\label{eq6.7}
\rho( cX) = c \rho (X) \text{ if } c  > 0 \text{ is a constant.}
\end{equation}
If \eqref{eq6.6} and \eqref{eq6.7} hold, then $\rho$ is called a {\it coherent} risk measure. See Artzner et al (1999) for details.

The case when $\alpha = \alpha_e$ is the entropic penalty case. It corresponds to the case when the driver is
$$g(z) ;= - \frac{1}{2} z^2$$
in the BSDE representation  \eqref{eq2.13} of convex risk measures, as explained in Theorem \ref{th2.13}. We will not prove this here, but we get an indirect confirmation of this connection by comparing Theorem \ref{th4.2} 
and the forthcoming Theorem \ref{th6.5} (based on \eqref{eq6.1} and $\alpha  = \alpha_e$).

\subsection{SD games and the HJBI equation}\label{sec6.2}

We now return to the problem of risk minimization in a financial market, but now the risk $\rho$ is represented by \eqref{eq6.1}. If $X_\varphi(t)$ is the wealth processes corresponding to a portfolio $\varphi$, then the risk minimization problem becomes
\begin{equation}\label{eq6.8}
\inf_{\varphi \in \AC} \left( \sup_{Q \in \PC} \{ E_Q [- X_\varphi(T)] - \alpha(Q)\} \right).
\end{equation}

This is a min-max problem which in our setting is called a {\it stochastic differential game}:

\begin{myenumerate}
\item If the system is Markovian, we can use an extension of the HJB equation to a corresponding Hamilton-Jacobi-Bellman-Isaacs (HJBI) equation for games.
\item We can extend the maximum principle approach to games. This approach does not require that the system is Markovian. See Remark  \ref{rem6.6}.
\end{myenumerate}

In this section we present the dynamic programming approach to stochastic differential games. We only present the the case with zero sum games. For extension to non-zero sum games, we refer to Mataramvura \& \O ksendal (2008).

Suppose the state $Y(t) = Y^{(u)}_y(t) \in \RB^k$ is described by a controlled SDE of the form
\begin{equation}\label{eq6.9}
\begin{cases}
dY(t) & \displaystyle = b(Y(t), u(t))dt + \sigma(Y(t), u(t)) dB(t) + \int_{\RB^l} \gamma(Y(t), u(t), \zeta) \tN (dt, d\zeta) \; ; \; t \geq 0 \\
Y(0) & = y \in \RB^k.
\end{cases}
\end{equation}

Here $u = (u_1, u_2)$ is a control process with values in $V_1 \times V_2 \subset \RB^{n_1} \times \RB^{n_2}$. The control $u_i$ is chosen by player number $i$, and $\AC_i$ is the set of admissible control processes for player number $i$; $i=1,2$. We put $\AC = \AC_1 \times \AC_2$. Let
$$f : \RB^k \times V_1 \times V_2 \rightarrow \RB$$
be a given profit rate and let
$$g : \RB^k \rightarrow \RB$$
be a given bequest function. We assume that
$$E^y \left[ \int_0^\infty |f(Y(t), u(t))| dt + |g(Y(\tau_{\Scal}))| \right] < \infty$$
for all $y$, where $E^y[\varphi(Y(t))]$ means $E[\varphi(Y_u(t))]$ etc...

Here
\begin{equation}\label{eq6.10}
\tau_{\Scal} := \inf \{ t > 0 \; ; \; Y^{(u)}_y(t) \not\in \Scal\},
\end{equation}
with $\Scal \in \RB^k$ being a given {\it solvency region}. We may interpret $\tau_{\Scal}$ as the {\it bankruptcy time}.

Let $\AC_1, \AC_2$ be given families of admissible control processes $u_1,u_2$ for player number 1 and 2, respectively

\begin{problem}[Zero-sum stochastic differential games]\label{prob6.1}

Find $\varphi(y)$ and $u^*_i \in \AC_{\epsilon_i} \; ; \; i=1,2$ such that
\begin{equation}\label{eq6.11}
\Phi(y) = \sup_{u_1 \in \AC_1} \left( \inf_{u_2 \in \AC_2} J^u(y)\right) = J^{u_1^*, u_2^*}(y),
\end{equation}
where
\begin{equation}\label{eq6.12}
J^{u}(y) = J^{u_1,u_2}(y) = E \left[ \int_0^{\tau_S} f (Y(s), u(s))ds + g(Y(\tau_S)) \chi_{\tau_S < \infty} \right]
\end{equation}
is the performance functional.
\end{problem}

As in classical Markovian stochastic control we can restrict ourselves to consider Markov controls (feedback controls), i.e. we assume that
$$u_i(t) = \tu_i(Y(t))$$
for some deterministic function $\tu_i : \RB^k \rightarrow \RB \; ; \; i=1,2$. As is customary we do not distinguish notationally between $u_i$ and $\tu_i$.

When the control $u = (u_1, u_2) \in \AC_1 \times \AC_2$ is Markovian, the corresponding controlled process $Y^{(u)}(t)$ will be a (Markovian) jump diffusion with generator $A^u$ given by
\begin{align}\label{eq6.13}
A^u \varphi(y) &= \sum_{i=1}^k b_i (y,u(y)) \frac{\partial \varphi}{\partial y_i} + \frac{1}{2} \sum_{i,j=1}^k (\sigma \sigma^T)_{ij} (y,u(y)) \frac{\partial^2 \varphi}{\partial y_i \partial y_j}(y) \nonumber \\
& + \sum_{j=1}^\ell \int_{\RB^k} \{ \varphi(y + \gamma^{(j)}(y,u(y), \zeta)) - \varphi(y) - \nabla \varphi(y) \gamma^{(j)}(y, u(y), \zeta)\} \nu_j(d\zeta)
\end{align}
for all smooth functions $\varphi : \RB^k \rightarrow \RB$ with compact support, where $\gamma^{(j)}$ denotes column number $j$ in the $k \times \ell$ matrix $\gamma = [\gamma_{ij}]_{\substack{1 \leq i \leq k \\ 1 \leq j \leq \ell}}$. 

For more information about stochastic control of jump diffusions we refer to \O ksendal \& Sulem (2007).

We now state the main result of this section.

\begin{theorem}[Mataramvura \& \O ksendal (2008)]\label{the6.2} (The Hamilton-Jacobi-Bellman-Isaacs (HJBI) equation).

Suppose there exists a function $\varphi \in C^2(\Scal) \times C(\RB^k)$ and a Markov control $\hu = (\hu_1, \hu_2) \in \AC_1 \times \AC_2$ such that
\begin{myenumerate}
\item $A^{\hu_1,v_2} \varphi(y) + f(y, \hu_1(y), v_2) \geq 0$ for all $v_2 \in V_2, \; y \in \Scal$
\item $A^{v_1,\hu_2} \varphi(y) + f(y, v_1, \hu_2(y)) \leq 0$ for all $v_1 \in V_1; \; y \in \Scal$
\item $A^{\hu_1, \hu_2} \varphi(y) + f(y, \hu_1(y), \hu_2(y)) = 0$ for all $y \in \Scal$
\item $\displaystyle \lim_{t \rightarrow \tau^-_{\Scal}} \varphi (Y^{(u)}(t)) = g(Y^{(u)}(\tau_{\Scal})) \chi_{\tau_S < \infty}$ a.e. for all $u \in \AC_1 \times \AC_2$,  $y \in \Scal$
\item The family $\{ \varphi(Y^{(u)}(\tau))\}_{\tau \in \TC}$ is unformly integrable, for all $y \in \Scal$, $u \in \AC_1 \times \AC_2$. Here $\TC$ denotes the set of all $\FC$-stopping times $\tau \leq \tau_{\Scal}$.
\end{myenumerate}
Then
\begin{align}\label{eq6.14}
\varphi(y) & = \Phi(y): = \sup_{u_1 \in \AC_1} \left( \inf_{u_2 \in \AC_2} J^{u}(y)\right) \nonumber \\
 & = \inf_{u_2 \in \AC_2} \left( \sup_{u_1 \in \AC_1} J^{u}(y)\right) = \sup_{u_1 \in \AC_1} J^{u_1, \hu_2}(y) \nonumber \\
 & = \inf_{u_2 \in \AC_2} J^{\hu_1, u_2}(y) = J^{\hu_1, \hu_2}(y)
 \end{align}
 and
 \begin{equation}\label{eq6.15}
 u^* := (\hu_1, \hu_2) \text{ is an optimal control.}
 \end{equation}
\end{theorem}

\dproof
Choose $u \in \AC$. Then by the Dynkin formula for jump diffusions (see e.g. \O ksendal \& Sulem (2007), Theorem 1.24) we have
$$ E^y[ \varphi(Y^u(\tau_{\Scal}^{(N)}))] = \varphi(y) + E^y \left[ \int_0^{\tau_\Scal^{(N)}} A^u\varphi(Y^u(t))dt\right]$$
where
$$\tau_\Scal^{(N)} = \tau_\Scal \wedge N \wedge \inf \{ t > 0 \; ; \; |Y(t)| \geq N\} \; ; \; N = 1,2, \ldots$$
Hence, by (i) we get, with $u_i(s) = u_i(Y(s))$,
$$E^y [ \varphi(Y^{\hu_1,u_2}(\tau_\Scal^{(N)} ))] \geq \varphi(y) - E^y \left[ \int_0^{\tau_S^{(N)} }f(Y^{\hu_1,u_2}(s), u_1(s), \hu_2(s))ds\right].$$
Therefore,
\begin{align}\label{eq6.16}
\varphi(y) &\leq E^y \left[ \int_0^{\tau_S^{(N)} } f(Y^{\hu_1,u_2}(s), \hu_1(s), u_2(s))ds + \varphi (Y^{\hu_1,u_2}(\tau_S^{(N)}))\right] \nonumber \\
& \rightarrow J^{\hu_1,u_2}(y) \text{ as } N \rightarrow \infty.
\end{align}
Here we have used condition (v) and the fact that $Y(\cdot)$ is quasi-left continuous (i.e. left continuous at stopping times; see Jacod \& Shiryaev (2003), Prop. I.2.26 and Prop. I.3.27). Since $u_2 \in \AC_2$ was arbitrary, we deduce from \eqref{eq6.16} that
\begin{equation} \label{eq6.17}
\varphi(y) \leq \inf_{u_2 \in \AC_2} J^{\hu_1,u_2}(y).
\end{equation}
It follows that
\begin{equation}\label{eq6.18}
\varphi(y) \leq \sup_{u_1 \in \AC_1} \left( \inf_{u_2 \in \AC_2} J^{u_1,u_2}(y)\right) = \Phi(y).
\end{equation}
Similarly, applying (ii) we get
\begin{equation}\label{eq6.19}
\varphi(y) \geq J^{u_1, \hu_2}(y) \text{ for all } u_1 \in \AC_1
\end{equation}
and therefore
\begin{equation}\label{eq6.20}
\varphi(y) \geq \sup_{u_1 \in \AC_1} J^{u_1,\hu_2}(y) \geq \inf_{u_2 \leq \AC_2} \left( \sup_{u_1 \in \AC_1} J^{u_1,u_2}(y)\right).
\end{equation}
In the same way, applying (iii) we get
\begin{equation}\label{eq6.21}
\varphi(y) = J^{\hu_1, \hu_2}(y).
\end{equation}
By combining \eqref{eq6.17}-\eqref{eq6.21} we obtain
\begin{align*}
\varphi(y) & = \inf_{u_2 \in \AC_2} J^{\hu_1,u_2}(y) \leq \sup_{u_1 \in \AC_1} \left( \inf_{u_2 \in \AC_2} J^{u_1,u_2}(y)\right) = \Phi(y) \\
& \leq \inf_{u_2 \in \AC_2} \left( \sup_{u_1 \in \AC_1} J^{u_1,u_2}(y)\right) \leq \sup_{u_1 \in \AC_1} J^{u_1, \hu_2}(y) \leq \varphi(y).
\end{align*}
We conclude that $\varphi(y) = \Phi(y)$, that \eqref{eq6.14} holds and that $u^* := (\hu_1, \hu_2)$ is optimal.
\fproof

\subsection{Entropic risk minimization by the HJBI equation}\label{sec6.3}

Suppose the financial market is as in \eqref{eq2.1}, but with $r=0$ and with jumps added, i.e.
\begin{equation}\label{eq6.22}
\begin{cases}
S_0(t) =1 \text{ for all } t \\
dS_1(t)  = \displaystyle S_1(t^-) \left[ \mu(t) dt + \sigma(t) dB(t) + \int_\RB \gamma(t,\zeta) \tN (dt,d\zeta) \right];  t  \geq 0\\
S_1(0) > 0.
\end{cases}
\end{equation}
If $\beta(t)$ is a self-financing portfolio representing the number of units of the risky asset (with unit price $S_1(t)$) held at time $t$, the corresponding wealth process $X(t) = X_\beta(t)$ will be given by
\begin{align}\label{eq6.23}
dX(t) & = \beta(t) dS_1(t) \nonumber \\
  & = w(t) \left[ \mu(t) dt + \sigma(t) dB(t) + \int_\RB \gamma(t,\zeta) \tN (dt,d\zeta)\right] \nonumber \\
  & = dX^w(t),
\end{align}
where $w(t) := \beta(t) S_1(t^-)$ is the {\it amount} held in the risky asset at time $t$.

If we use the representation \eqref{eq6.1} of the risk measure $\rho = \rho_e$ corresponding to the family $\PC_\Theta$ of measures $Q$ given in \eqref{eq6.2}-\eqref{eq6.4} and the entropic penalty $\alpha_e$ given by \eqref{eq6.5}, the risk minimizing portfolio problem \eqref{eq6.8} can be written
\begin{equation}\label{eq6.24}
\inf_{w \in \WC} \left( \sup_{\theta \in \Theta} \left\{E \left[ - \frac{dQ}{dP} X^w(T) - \frac{dQ}{dP} log \left( \frac{dQ}{dP}\right) \right]\right\} \right)
\end{equation}
where $\WC$ is the family of admissible portfolios $w$. To put this problem into the setting of Section \ref{sec6.2}, we represent $Q$ by $\frac{dQ}{dP} = M^\theta(T)$ and we put
\begin{align}\label{eq6.25}
dY(t) & = dY^{\theta,w}(t) = \begin{bmatrix} dt \\ dX^w (t) \\ dM^\theta(t)\end{bmatrix} \nonumber \\
& = \begin{bmatrix} 1 \\ w(t) \mu(t) \\ 0\end{bmatrix} dt + \begin{bmatrix} 0 \\ w(t) \sigma(t) \\ M^\theta(t) \theta_0(t)\end{bmatrix} dB(t) + \int_\RB \begin{bmatrix} 0 \\ w(t) \gamma(t, \zeta) \\ M^\theta(t^-) \theta_1(t_1,\zeta)\end{bmatrix} \tN(dt, d \zeta),
\end{align}
with initial value
\begin{equation}\label{eq6.26}
Y^{\theta,w}(0) = y = \begin{bmatrix} s \\ x \\ m\end{bmatrix} \; ; \; s \in [0,T], x > 0, m > 0.
\end{equation}
In this case the solvency region is $\Scal = [0,T] \times \RB_+ \times \RB_+$ and the performance functional is
\begin{equation}\label{eq6.27}
J^{\theta,w}(s,x,m) = E^{s,x,m} [ - M^\theta(T) X^w(T) - M^\theta(T) \log M^\theta(T)].
\end{equation}
Assume from now on that
\begin{equation}\label{eq6.28}
\mu(t), \sigma(t) \text{ and } \gamma (t,\zeta) \text{ are deterministic.}
\end{equation}
Then $Y^{\theta,w}(t)$ becomes a controlled jump diffusion, and the risk minimization problem \eqref{eq6.8} is the following special case of Problem \ref{prob6.1}:
\begin{problem}[Entropic risk minimization]\label{prob6.3}

Find $w^* \in \WC$, $\theta^* \in \Theta$ and $\Phi(y)$ such that
\begin{equation}\label{eq6.29}
\Phi(y) = \inf_{w \in \WC} \left( \sup_{\theta \in \Theta} J^{\theta,w}(y) \right) = J^{ \theta^*, w^*}(y) \; ; \; y \in \Scal.
\end{equation}
\end{problem}
By \eqref{eq6.25} we see that the generator $A^{\theta,w}$ is given by (see \eqref{eq6.13})
\begin{align}\label{eq6.30}
A^{\theta,w} &\varphi (s,x,m) = \frac{\partial \varphi}{\partial s} (s,x,m) + w \mu(s) \frac{\partial \varphi}{\partial x}(s,x,m) \nonumber \\
& + \frac{1}{2} w^2 \sigma^2(s) \frac{\partial^2 \varphi}{\partial x^2}(s,x,m) + \frac{1}{2} m^2 \theta^2_0 \frac{\partial^2 \varphi}{\partial m^2}(s,x,m) + w \theta_0 m \sigma(s) \frac{\partial^2 \varphi}{\partial x \partial m} (s,x,m) \nonumber \\
 & + \int_\RB \{ \varphi(s,x+w \gamma(s,\zeta), m + m \theta_1(\zeta)) - \varphi(s,x,m) \nonumber \\
  &\left. - \frac{\partial \varphi}{\partial x} (s,x,m) w \gamma (s,\zeta) - \frac{\partial \varphi}{\partial m} (s,x,m) m \theta_1 (\zeta) \right\} \nu (d \zeta).
\end{align}
Comparing with the general formulation in Section \ref{sec6.3}, we see that in this case
$$f = 0 \text { and }  g(y) = g(x,m) = - m x - m \log (m).$$
Therefore, according to Theorem \ref{the6.2}, we should try to find a function $\varphi(s,x,m) = \rho^2(\Scal) \cap C(\RB^3)$ and control values $\theta = \hth(y)$, $w = \hw(y)$ such that
\begin{equation}\label{eq6.31}
\inf_{w \in \RB} \left( \sup_{\theta \in \RB^2} A^{\theta,w} \varphi(y)\right) = A^{\hth, \hw} \varphi(y) \; ; \; y \in \Scal
\end{equation}
and
\begin{equation}\label{eq6.32}
\lim_{t \rightarrow T^-} \varphi(s,x,m) = - xm - m \log(m).
\end{equation}
Let us try a function of the form
\begin{equation}\label{eq6.33}
\varphi(s,x,m) = - xm - m \log(m) + \kappa(s)m
\end{equation}
where $\kappa$ is a deterministic function, $\kappa(T) = 0$. Then by \eqref{eq6.30}
\begin{align}\label{eq6.34}
A^{\theta,w} & \varphi(s,x,m) = \kappa'(s)m + m \mu(t) (-m) + \frac{1}{2} m^2 \theta^2_0 \left( - \frac{1}{m}\right) + w \theta_0 m \sigma(t) (-1) \nonumber \\
& \quad +  \int_\RB \{ - (x + w \gamma(s,\zeta))(m + m \theta_1(\zeta)) + xm - (m + m \theta_1(\zeta)) \log(m + m \theta_1(\zeta)) \nonumber \\
& \quad + m \log m  + \kappa(s) (m + m \theta_1(\zeta)) - \kappa(s) m + mw \gamma(s,\zeta) \nonumber \\
& \quad - m \theta_1(\zeta) (- x - 1 - \log m + \kappa(s))\} \nu (d\zeta) \nonumber \\
& = m \left[\kappa'(s) - w \mu (t) - \frac{1}{2} \theta^2_0 - w \theta_0 \sigma(t)
 + \int_\RB \theta_1(\zeta) \{ 1 - \log (1 + \theta_1(\zeta)) - w \gamma (s,\zeta)\} \nu (d\zeta) \right]
\end{align}
Maximizing $A^{\theta,w}\varphi(y)$ with respect to $\theta = (\theta_0, \theta_1)$ and minimizing with respect to $w$ gives the following first order equations
\begin{equation}\label{eq6.35}
\hth_0(s) + \hw(s) \sigma(s) = 0
\end{equation}
\begin{equation}\label{eq6.36}
1 - \log (1 + \hth_1(s,\zeta)) - \hw(s) \gamma(s,\zeta) - \frac{\hth_1(s,\zeta)}{1 + \hth_1(s,\zeta)} = 0
\end{equation}
\begin{equation}\label{eq6.37}
\mu(s) + \hth_0(s) \sigma(s) - \int_\RB \hth_1(s,\zeta) \gamma(s,\zeta) \nu (d \zeta) = 0.
\end{equation}
These are 3 equations in the 3 unknown candidates $\hth_0, \hth_1$ and $\hw$ for the optimal control for the SD game in Problem \ref{prob6.3}.
To get an explicit solution, let us now assume that
\begin{equation}\label{eq6.38}
N=0 \text{ and } \gamma = \theta_1 = 0.
\end{equation}
Then \eqref{eq6.35}-\eqref{eq6.37} gives
\begin{equation}\label{eq6.39}
\hth_0(s) = - \frac{\mu(s)}{\sigma(s)}, \hw(s) = \frac{\mu(s)}{\sigma^2(s)}.
\end{equation}
Substituted into \eqref{eq6.34} we get by (iii) Theorem \ref{the6.2}
$$A^{\hth, \hw} \varphi(s,x,w) = m \left[ \kappa'(s) - \frac{1}{2} \left( \frac{\mu(s)}{\sigma(s)}\right)^2\right] = 0.$$
Combining this with (iv) Theorem \ref{the6.2} we obtain
\begin{equation}\label{eq6.40}
\kappa(s) = - \int_s^T \frac{1}{2}\left( \frac{\mu(t)}{\sigma(t)}\right)^2 dt.
\end{equation}
Now all the conditions of Theorem \ref{the6.2} are satisfied, and we get:
\begin{theorem}[Entropic risk minimization]\label{th6.5}

Assume that \eqref{eq6.28} and \eqref{eq6.38} hold. Then the solution of Problem \ref{prob6.3} is
\begin{equation}\label{eq6.41}
\Phi(s,x,m) = - xm - m \log m - \int_s^T \frac{1}{2} \left( \frac{\mu(t)}{\sigma(t)}\right)^2dt
\end{equation}
and the optimal controls are
\begin{equation}\label{eq6.42}
\hth_0(s) = -  \frac{\mu(s)}{\sigma(s)} \text{ and } \hw(s) = \frac{\mu(s)}{\sigma^2(s)}  \; ; \; s \in [0,T].
\end{equation}
In particular, choosing the initial values $s=0$ and $m=1$ we get
\begin{equation}\label{eq6.43}
\Phi (0,x,1) = - x - \int_0^T \frac{1}{2}  \left( \frac{\mu(t)}{\sigma(t)}\right)^2 dt.
\end{equation}
This agrees with what we obtained in Theorem \ref{th4.2}, using the maximum principle for FBSDES.   
\end{theorem} 
\begin{remark}\label{rem6.6}
\begin{itemize}
\item

Recently, in \O ksendal \& Sulem (2012), a maximum principle for stochastic differential games has been developed, along the same lines as in Section \ref{sec3}. This principle could have been used as an alternative approach to risk minimization when the risk is given by the dual representation in Section \ref{sec6.1}. 
\item
Theorem \ref{th6.5} could also have been obtained by using a stochastic HJB approach to the optimal control of forward-backward SDEs. See Theorem 3.4 in \cite{OSZ1}.
\end{itemize}
\end{remark}

\paragraph{Acknowledgements}
These lecture notes are based on our earlier works and lectures given on this topic, and on a course on risk minimization that B.\O . gave at NHH,  in  2013. We are grateful to K. Aase, J.  Haug, S.-A. Persson and J. Ub\o e for valuable comments.

 \end{document}